# $(k, s)$-POSITIVITY AND VANISHING THEOREMS FOR COMPACT KAHLER MANIFOLDS [*]

Qilin Yang [†][‡]

**Abstract**

We study the $(k, s)$-positivity for holomorphic vector bundles on compact complex manifolds. $(0, s)$-positivity is exactly the Demailly $s$-positivity and a $(k, 1)$-positive line bundle is just a $k$-positive line bundle in the sense of Sommese. In this way we get a unified theory for all kinds of positivities used for semipositive vector bundles. Several new vanishing theorems for $(k, s)$-positive vector bundles are proved and the vanishing theorems for $k$-ample vector bundles on projective algebraic manifolds are generalized to $k$-positive vector bundles on compact Kähler manifolds.

## 1 Introduction

Vanishing theorems play an important role in the study of algebraic and analytic varieties by cohomological methods. The celebrated Kodaira vanishing theorem was proved originally by using the Bochner techniques in differential geometry [21, 46]. Using his vanishing theorem, Kodaira proved the projective embedding theorem now bearing his name, which says that a compact Kähler manifold is projective algebraic if and only if it admits an ample, or equivalently, a positive line bundle. Followed by Chow's theorem, the Kodaira vanishing theorem sets up a bridge between the analytic geometry and the algebraic geometry. The Kodaira vanishing theorem was further studied in many ways. Generalizations have been found, on the differential geometric side, by Akizuki-Nakano [1], Griffiths [17], Girbau [15], Gigante [14], Ohsawa-Takegoshi [33, 34, 35, 36, 43] and Demailly [7, 8] et al; on the algebraic geometric side, by Mumford [30], Grauert-Riemenschneider [16], Sommese [41], Kawamata [19], Viehweg [44] and Lazarsfeld [23] et al. We refer readers to [12, 23] for a survey of history and developments and applications by algebraic approach.

In this paper we will investigate vanishing theorems from the viewpoint of complex differential geometry. Our original motivation was to understand the Demailly $s$-positivity (cf. [7]), which is equivalent to the Griffiths positivity when $s = 1$ and to the Nakano positivity when $s$ attains its maximum. The Demailly $s$-positivity is a concept which is difficult to understand and to check. In Section 3, we introduce the notion of $(k, s)$-positivity from the viewpoint of the theory of hermitian forms in linear algebra. When $k = 0$, it is equivalent to Demailly's $s$-positivity; if $s = 1$ and in the case of line bundle, it is just $k$-positivity used implicitly by Gigante and Girbau in [14, 15]. We present some properties of $(k, s)$-positivity and then give some vanishing theorems for $(k, s)$-positive vector bundles (cf. Theorem 3.9 and Corollary 3.13) which generalize those obtained by Demailly in [7, 9] for $s$-positive vector bundles in his sense.

In 1978, Sommese introduced the concept of $k$-ampleness for vector bundles on any complex space, in the case of line bundle the 0-ampleness being ampleness in the usual sense when the base space is projective. He established several vanishing theorems for $k$-ample vector bundles on projective algebraic manifolds [41]. Surprisingly, the same results are still hold for $k$-positive vector bundles on compact Kähler manifolds, (cf.[40] and [20]). In fact, in a later paper [42], Sommese proved that $k$-ample vector bundles are $k$-positive, but the statement was expressed in another way and not easy to follow. For the reader's convenience, we present a sketch of the proof that $k$-ampleness implies $k$-positivity in Section 4. In Section 5, we generalize Theorem 3.9 to the case of flag bundles (cf. Theorem 5.3) on compact complex manifolds, which says that the product of tautological quotient determinant bundles of the flag bundle of a holomorphic vector bundle $E$ is $k$-positive if $E$ itself is Griffiths $k$-positive. Using Theorem 5.3 we generalize the Demailly vanishing theorems for tensor powers of ample bundles on projective algebraic manifolds to the Griffiths $k$-positive vector bundles on compact Kähler manifolds (cf. Theorems 5.10, 5.11). Finally, in Section 6, we present some questions motivated by the results we get form Section 2 to Section 5.

## 2 Bochner-Kodaira-Nakano inequality

In this section we will recall the famous Bochner-Kodaira-Nakano inequality, which is the starting point for all kinds of positivity that leads to the vanishing theorems. We will present the Nakano vanishing theorem and give a simplified proof of the Gigante-Girbau vanishing theorem since it not only includes the main ingredients in using the Bochner-Kodaira-Nakano inequality but also will be used again and again in this paper.

Let $X$ be a complex manifold of dimension $n$ and $E$ a holomorphic vector bundle of rank $r$ over $X$. The Dolbeault operator $\bar{\partial}$ acts linearly on the space $\Omega_X^{p,q}(E)$ of smooth $E$-valued $(p,q)$-forms; if $f$ is a smooth (p,q)-form and $s$ is a holomorphic section of $E$, then define $\bar{\partial}(f \otimes s) = \bar{\partial}f \otimes s$. The Dolbeault cohomology groups $H^{p,q}(X, E)$ are defined to be the cohomology groups of the complex

$$\cdots \xrightarrow{\bar{\partial}} \Omega_X^{p,q}(E) \xrightarrow{\bar{\partial}} \Omega_X^{p,q+1}(E) \xrightarrow{\bar{\partial}} \cdots.$$

The soft sheaves $\Omega_X^{p,\bullet}(E)$ give an acyclic resolution of the locally free sheaf $\Omega_X^p(E)$ of holomorphic $E$-valued $p$-forms. Hence we have the following canonical isomorphism (called the *Dolbeault isomorphism*):

$$H^{p,q}(X,E) \cong H^q(X, \Omega_X^p(E)).$$

Now suppose that $X$ is equipped with a hermitian metric $g$, and $E \xrightarrow{\pi} X$ is any complex vector bundle with a hermitian metric $h$ and a hermitian connection $D$, which splits in a unique way as a sum of a $(1,0)$ and a $(0,1)$ connection: $D = D' + D''$. If $E$ is a holomorphic vector bundle, there is a unique hermitian connection, called *Chern* connection such that $D'' = \bar{\partial}$ and hence $D'^2 = D''^2 = 0$. Thus the curvature form $\Theta_h(E) = D'D'' + D''D'$ is a $End(E)$-valued $(1,1)$-form, and $i\Theta_h(E)$ is called the *Chern Curvature form* of $E$. Let $\{z^j\}$ be the local holomorphic coordinate of $X$, $\{e_\alpha\}$ a holomorphic orthonormal frame and $\{e^\alpha\}$ the dual orthonormal frame of $E$. Let $(g_{j\bar{k}})$ and $(h_{\alpha\bar{\beta}})$ be the hermitian metrics on $X$ and on $E$ respectively, and their inverses denoted respectively by $(g^{i\bar{j}})$ and $(h^{\alpha\bar{\beta}})$. Then we can write the Kähler form and the curvature form respectively as

$$\omega = ig_{j\bar{k}}dz^j \wedge d\bar{z}^k, \quad i\Theta_h(E) = R^\alpha_{\beta j\bar{k}}dz^j \wedge d\bar{z}^k \otimes e_\alpha \otimes e^\beta, \tag{2.1}$$

where

$$R^\alpha_{\beta j\bar{k}} = -h^{\alpha\bar{\gamma}}\partial_j\bar{\partial}_k h_{\beta\bar{\gamma}} + h^{\alpha\bar{\gamma}}h^{\lambda\bar{\mu}}\partial_j h_{\beta\bar{\mu}}\bar{\partial}_k h_{\lambda\bar{\gamma}}. \tag{2.2}$$

Here and throughout the rest of the paper we use the summation convention of summing over any index which appears once as a subscript and once as a superscript if without special mentions. The first Chern class $c_1(E) \in H^2(X, \mathbb{R})$ could be represented by the first Chern curvature form

$$c_1(E, h) = \text{Tr}_E(\frac{i}{2\pi}\Theta_h(E)) = \frac{i}{2\pi}\Theta_h(\det(E)). \tag{2.3}$$

Conversely, if $X$ is a compact Kähler manifold, any 2-form representing the first Chern class $c_1(E)$ is in fact the first Chern curvature form $c_1(E, h)$ for some hermitian metric $h$ on $E$ (cf.[20]). In local coordinate we have

$$i\Theta_h(\det(E)) = R_{j\bar{k}}dz^j \wedge d\bar{z}^k = -\partial\bar{\partial}\log\det(h_{\alpha\bar{\beta}}) \tag{2.4}$$

with $R_{j\bar{k}} = R^\alpha_{\alpha j\bar{k}}$. Note that if $E$ is a line bundle then its curvature form represents its first Chern class up to multiplying by $\frac{1}{2\pi}$.

Let $u = u^\alpha_{J_p, \bar{K}_q}dz^{J_p} \wedge d\bar{z}^{\bar{K}_q} \otimes e_\alpha$ and $v = v^\alpha_{J_p, \bar{K}_q}dz^{J_p} \wedge d\bar{z}^{\bar{K}_q} \otimes e_\alpha$ be $E$-valued $(p,q)$-forms, where $J_p = (j_1, \cdots, j_p)$ and $\bar{K}_q = (\bar{k}_1, \cdots, \bar{k}_q)$ respectively run all $p$-tuples and $q$-tuples, and $dz^{J_p} = dz^{j_1} \wedge \cdots \wedge dz^{j_p}$ and $d\bar{z}^{\bar{K}_q} = d\bar{z}^{k_1} \wedge \cdots \wedge d\bar{z}^{k_q}$. Define their point-wise inner product by

$$\langle u, v \rangle_\omega = \frac{1}{p!q!}g^{J_p\bar{L}_p}g^{M_q\bar{K}_q}u^\alpha_{J_p\bar{K}_q}h_{\alpha\bar{\beta}}\overline{v^\beta_{L_p\bar{M}_q}} = \frac{1}{p!q!}u^\alpha_{J_p\bar{K}_q}h_{\alpha\bar{\beta}}\overline{v^{\beta J_p\bar{K}_q}}, \tag{2.5}$$

and $|u|_\omega^2 = \langle u, u \rangle_\omega$ the corresponding norm square of $u$, where $g^{J_p \bar{L}_p} = g^{j_1 \bar{l}_1} \cdots g^{j_p \bar{l}_p}$. The $L^2$ inner product of $u$ and $v$ is defined by

$$(u, v)_\omega = \int_X \langle u, v \rangle_\omega dV_\omega, \tag{2.6}$$

where $dV_\omega = \frac{\omega^n}{n!}$. The induced $L^2$ norm on $u$ is denoted by $||u||_\omega^2 = (u, u)_\omega$.

If $X$ is compact, then complex Laplace operator $\Delta'' = D'' D''^* + D''^* D'' : \Omega_X^{p,q}(E) \to \Omega_X^{p,q}(E)$, is a self-adjoint elliptic operator. Here $D''^*$ is the adjoint of $D''$ with respect to the $L^2$ inner product. The Hodge theorem states that the space of harmonic $(p,q)$-forms with values in $E$, i.e.,

$$H^{p,q}(E) = \mathrm{Ker}(\Delta'' : \Omega_X^{p,q}(E) \to \Omega_X^{p,q}(E)),$$

is finite dimensional, and there are isomorphisms

$$H^{p,q}(E) \cong H^{p,q}(X, E).$$

Let us now suppose that $X$ is a compact Kähler manifold, i.e., the Kähler form $\omega$ of $X$ is closed: $d\omega = 0$. Let $L : \Omega_X^{p,q}(E) \to \Omega_X^{p+1,q+1}(E)$ be the operator defined by $Lu = \omega \wedge u$ and $\Lambda = L^*$ is the adjoint operator. In this case the complex Laplace operators $\Delta' = D'D'^* + D'^*D'$ and $\Delta'' = D''D''^* + D''^*D''$ are related by the following Bochner-Kodaira-Nakano identity

$$\Delta'' = \Delta' + [i\Theta_h(E), \Lambda]. \tag{2.7}$$

If $u$ is an arbitrary $E$-valued $(p,q)$-form, then an integration by part yields $(\Delta'u, u)_\omega = ||D'u||_\omega^2 + ||D'^*u||_\omega^2$ and $(\Delta''u, u)_\omega = ||D''u||_\omega^2 + ||D''^*u||_\omega^2$, therefore we get the following inequality

$$||D''u||_\omega^2 + ||D''^*u||_\omega^2 \geq \int_X \langle [i\Theta_h(E), \Lambda]u, u \rangle_\omega dV_\omega. \tag{2.8}$$

This inequality is known as the *Bochner-Kodaira-Nakano* inequality. If $u \in H^{p,q}(E)$, then $u$ is $\Delta''$-harmonic and hence $D''u = D''^*u = 0$ by the Hodge theory. Furthermore if that $[i\Theta_h(E), \Lambda]$ is positive everywhere on $\Omega^{p,q}(E)$, then $u = 0$. Hence $H^{p,q}(X, E) \cong H^q(X, \Omega_X^p(E)) \cong H^{p,q}(E) = 0$. Thus we get a vanishing cohomology group. Therefore to prove a vanishing theorem for $E$-valued Dolbeault cohomology group, the key point is to find conditions under which the operator $[i\Theta_h(E), \Lambda]$ is positive definite.

For $u = u_{J,\bar{K}}^\alpha dz^J \wedge dz^{\bar{K}} \otimes e_\alpha \in \Omega_X^{p,q}(E)$, we have the following formula:

$$\begin{aligned}
\langle [i\Theta_h(E), \Lambda]u, u \rangle_\omega &= \tfrac{1}{(p-1)!q!} R_{\alpha\bar{\beta}j\bar{k}} g^{l\bar{k}} u_{lR_{p-1},\bar{K}_q}^\alpha \overline{u^{\beta,\bar{j}\bar{R}_{p-1},K_q}} \\
&+ \tfrac{1}{p!(q-1)!} R_{\alpha\bar{\beta}j\bar{k}} g^{j\bar{l}} u_{J_p,\bar{l}\bar{S}_{q-1}}^\alpha \overline{u^{\beta,\bar{J}_p,kS_{q-1}}} \\
&- \tfrac{1}{p!q!} R_{\alpha\bar{\beta}j\bar{k}} g^{j\bar{k}} u_{J_p,\bar{K}_q}^\alpha \overline{u_{\bar{J}_p,K_q}^\beta}.
\end{aligned} \tag{2.9}$$

If we choose normal metrics both on the base $X$ and the fibre $E$ with $g_{j\bar{k}} = \delta_{j\bar{k}}$ and $h_{\alpha\bar{\beta}} = \delta_{\alpha\bar{\beta}}$, then we have the following simpler expression:

$$\begin{aligned}\langle [i\Theta_h(E), \Lambda]u, u\rangle_\omega &= \sum R_{\alpha\bar{\beta}j\bar{k}} u^\alpha_{kR_{p-1},\bar{K}_q} \overline{u^\beta_{jR_{p-1},K_q}} \\ &+ \sum R_{\alpha\bar{\beta}j\bar{k}} u^\alpha_{J_p,\bar{j}\bar{S}_{q-1}} \overline{u^\beta_{J_p,kS_{q-1}}} \\ &- \sum R_{\alpha\bar{\beta}j\bar{j}} u^\alpha_{J_p,\bar{K}_q} \overline{u^\beta_{J_p,K_q}},\end{aligned} \qquad (2.10)$$

where the summation is on all indices $1 \leq j, k \leq n, 1 \leq \alpha, \beta \leq r$ and on all multi-indices $J_p, K_q, R_{p-1}, S_{q-1}$ of an increasing order with $|J_p| = p, |K_q| = q, |R_{p-1}| = p-1$ and $|S_{q-1}| = q - 1$. However, it is still very hard to decide when the expression (2.10) is positive, except in some special cases.

A first tractable case is that when $p = n$. Then in the first summation of (2.10) we must have $j = k$ and $R_{n-1} = \{1, \cdots, n\} \setminus \{j\}$. So the first summation cancels the last summation.

**Definition 2.1.** Let $E$ be a holomorphic vector bundle on a compact complex manifold $X$. $E$ is called *Nakano positive*, if

$$i\Theta_h(E)(u, u) = R_{\alpha\bar{\beta}j\bar{k}} u^{\alpha j} \overline{u^{\beta k}} > 0, \qquad (2.11)$$

for any non-zero tensor $u = u^{j\alpha}\partial/\partial z_j \otimes e_\alpha \in TX \otimes E$; $E$ is called *Griffiths positive* if

$$i\Theta_h(E)(u, u) = R_{\alpha\bar{\beta}j\bar{k}} \xi^j \bar{\xi}^k v^\alpha \bar{v}^\beta > 0 \qquad (2.12)$$

for all non-zero decomposable tensors $u = \xi \otimes v \in TX \otimes E$.

If the inequalities in (2.11) and (2.12) are not strict, $E$ is called *Nakano semipositive* and *Griffiths semipositive* respectively. These definitions were introduced by Kodaira [21], Nakano [32] and Griffiths [17]. It is clear that Griffiths positivity is weaker than Nakano positivity. Though Griffiths positivity does not directly produce vanishing theorems, it has nicer functorial properties and is very close to the concept of ampleness for algebraic vector bundles. We will discuss their relations in the next section.

If $E$ is Nakano-positive, then $[i\Theta_h(E), \Lambda]$ is positive definite on $E$-valued $(n, q)$-forms for any $q \geq 1$ and we have (cf. [32]):

**Theorem 2.2.** (Nakano's vanishing theorem). *If $E$ is a Nakano positive holomorphic vector bundle on a compact Kähler manifold $X$, then we have*

$$H^{n,q}(X, E) = H^q(X, K_X \otimes E) = 0, \quad q \geq 1.$$

Another tractable case is when $E$ is a holomorphic line bundle, in this case we denote $E$ by $B$. Then we can choose a coordinate system at each point $x \in X$, which diagonalizes simultaneously the hermitian form $\omega(x)$ and $i\Theta_h(B)|_x$, such that

$$\omega(x) = i\sum_{j=1}^n \mu_j(x) dz^j \wedge d\bar{z}^j, \qquad i\Theta_h(B)|_x = i\sum_{j=1}^n \nu_j(x) dz^j \wedge d\bar{z}^j.$$

Without loss of generality, assume that $\nu_1(x)/\mu_1(x) \leq \cdots \leq \nu_n(x)/\mu_n(x)$. Then for any $(p,q)$-form $u = u_{J,\bar{K}} dz^J \wedge z^{\bar{K}} \otimes e$, from (2.10) we have

$$
\begin{aligned}
& \langle [i\Theta_h(B), \Lambda] u, u \rangle_\omega(x) \\
= & \sum_{|J|=p, |K|=q} \left( \sum_{j \in J} \frac{\nu_j(x)}{\mu_j(x)} + \sum_{j \in K} \frac{\nu_j(x)}{\mu_j(x)} - \sum_{j=1}^n \frac{\nu_j(x)}{\mu_j(x)} \right) |u_{J,\bar{K}}|_\omega^2 \\
\geq & \left( \frac{\nu_1(x)}{\mu_1(x)} + \cdots + \frac{\nu_p(x)}{\mu_p(x)} + \frac{\nu_1(x)}{\mu_1(x)} + \cdots + \frac{\nu_q(x)}{\mu_q(x)} - \sum_{j=1}^n \frac{\nu_j(x)}{\mu_j(x)} \right) |u|_\omega^2.
\end{aligned} \quad (2.13)
$$

Observe that if the ratio $\nu_j(x)/\mu_j(x)$ varies small when $j$ varies, for example, in the extreme case when all $\nu_j(x)/\mu_j(x)$ are equal, then $[i\Theta_h(E), \Lambda](x)$ is positive when $p+q > n$. This observation tells us that if we choose the Kähler metric $\omega$ properly such that the eigenvalues of $i\Theta_h(E)$ vary mildly relative to $\omega$, then we can deduce the positivity of $[i\Theta_h(E), \Lambda]$.

**Definition 2.3.** (cf. [20, 40]) A holomorphic line bundle $B$ on a compact complex manifold $X$ is called $k$-*positive*, if there is a hermitian metric on $B$ such that its first Chern curvature form $c_1(B, h)$ is semi-positive and has at least $n - k$ positive eigenvalues.

Since the first Chern curvature form coincides with its curvature form up to a positive constant by (2.3), $B$ is $k$-positive if and only if $i\Theta_h(B)$ is semi-positive with rank $i\Theta_h(B) \geq n - k$. In this case we can make a special choice of the Kähler metric $\omega_\kappa = i\Theta_h(B) + \kappa\omega$ for some positive number $\kappa$. Let $\Lambda_\kappa$ be the associated adjoint operator of multiplication by $\omega_\kappa$. Then the eigenvalues of $i\Theta_h(B)$ relative to $\omega_\kappa$ are $\{r_j(x)\}_{1 \leq j \leq n}$ with

$$
r_j(x) = \frac{\nu_j(x)}{\kappa \mu_j(x) + \nu_j(x)} = \frac{\nu_j(x)/\mu_j(x)}{\kappa + \nu_j(x)/\mu_j(x)} = \frac{\frac{\nu_j(x)/\mu_j(x)}{\kappa}}{1 + \frac{\nu_j(x)/\mu_j(x)}{\kappa}} = 1 - \frac{1}{1 + \frac{\nu_j(x)/\mu_j(x)}{\kappa}}.
$$

Fix a point $x \in X$ and assume that $\text{rank}(\Theta_h(B)|_x) = n - s \geq n - k$. Then $0 = r_s(x) < r_{s+1}(x) \leq \cdots \leq r_n(x)$. Thus $r_j(x) = 0$ for $j \leq s$. If we choose $\kappa \to 0^+$ then $r_j(x) \to 1$ for all $s+1 \leq j \leq n$. If $p + q > n + k$, then $s \leq \min\{p, q\}$ and

$$
\begin{aligned}
& \lim_{\kappa \to 0^+} (r_1(x) + \cdots + r_p(x) + r_1(x) + \cdots + r_q(x) - (r_1(x) + \cdots + r_n(x))) \\
= & \ p - s + q - s - (n - s) \\
= & \ (p + q) - (n + s) \geq (p + q) - (n + k) \geq 1.
\end{aligned}
$$

Since $X$ is compact, we can use a finite cover by open sets, such that on each open set, if $\kappa$ is a sufficiently small positive number we may have on each open neighborhood and hence everywhere on $X$ that

$$
\langle [i\Theta_h(B), \Lambda_\kappa] u, u \rangle_{\omega_\kappa} \geq \frac{1}{2} |u|_{\omega_\kappa}^2.
$$

As a consequence we have (cf. [14, 15]):

**Theorem 2.4.** (Gigante-Girbau's vanishing theorem). *If $B$ is a $k$-positive line bundle on a compact Kähler manifold $X$, then*

$$
H^{p,q}(X, B) = H^q(X, \Omega_X^p \otimes B) = 0, \quad \text{for} \quad p + q > n + k.
$$

If $B$ is a positive line bundle, the above theorem was firstly proved by Akizuki-Nakano in [1]. This theorem is further generalized by Ohsawa-Takegoshi in [43] to weakly 1-complete Kähler manifold.

**Remark 2.5.** (1). Let $X = \mathbb{P}^k \times \mathbb{P}^{n-k}$ and $B = \pi_2^* \mathcal{O}_{\mathbb{P}^{n-k}}(1)^{\otimes N}$, where $\pi_2$ is the projection from $\mathbb{P}^k \times \mathbb{P}^{n-k}$ to its second factor, $N$ is a sufficiently large positive integer and $\mathcal{O}_{\mathbb{P}^{n-k}}(1)$ is the hyperplane bundle of $\mathbb{P}^{n-k}$. Then $B$ is $k$-positive but is not $k-1$-positive on $X$. We have

$$0 \neq H^k(\mathbb{P}^k, \Omega_{\mathbb{P}^k}^k) \otimes H^0(\mathbb{P}^{n-k}, \Omega_{\mathbb{P}^{n-k}}^{n-k}(\mathcal{O}_{\mathbb{P}^{n-k}}(N))) \subset \bigoplus_{p+q=n+k} H^p(X, \Omega_X^q(B)).$$

Thus the Nakano and Gigante-Girbau vanishing theorems are optimal. The following example (due to Ramanujam (cf. [40], pp. 55) shows that in the Theorem 2.4, it is not enough to assume that $B$ is $k$ positive at an open set of $X$ instead of at all points of $X$: For $n \geq 3$ let $\pi : \widetilde{\mathbb{P}^n} \mapsto \mathbb{P}^n$ be the blowing up at a point. Then the pull back of the hyperplane bundle $\pi^* \mathcal{O}_{\mathbb{P}^n}(1)$ is a semipositive line bundle, it has maximal rank on an open set of $\widetilde{\mathbb{P}^n}$. However $H^{n-1}(\widetilde{\mathbb{P}^n}, \Omega_{\widetilde{\mathbb{P}^n}}^{n-1} \otimes \pi^* \mathcal{O}_{\mathbb{P}^n}(1)) \neq 0$. Combining the Ramanujam's example with similar product construction used above, we can also give examples where the line bundle is semipositive and $k$-positive on an open set for any $k > 0$.

(2). Recall that a line bundle $B$ is called *nef* if there exists a smooth hermitian metric $h_\epsilon$ such that $i\Theta_{h_\epsilon}(B) \geq -\epsilon\omega$ for any $\epsilon > 0$. A semipositive line bundle is nef but the converse is not true. However, a nef line bundle is, in an approximate sense, very close to a semipositive line bundle. Could we get some vanishing Dolbeault cohomology groups twisted by a nef line bundle? In the proof of Theorem 2.4, the hermitian form $i\Theta_h(B)$ plays a dominant role in the sequence of Kähler metrics $\omega_\kappa = i\Theta_h(B) + \kappa\omega$ as $\kappa \to 0$. To use the trick of the proof of Theorem 2.4, the line bundle $B$ must be semipositive at least. Thus to answer the question we should look for new ideas.

(3). The celebrated Kawamata-Viehweg vanishing theorem has been generalized to Kähler manifolds by Nadel, which states that

$$H^p(X, K_X \otimes B \otimes \mathscr{I}(h)) = 0,$$

for $p \geq 1$ for a pseudo-effective line bundle $B$ on a Kähler manifold with singular hermitian metric $h$ whose curvature form as a current is strictly positive, where $\mathscr{I}(h) \subset \mathcal{O}_X$ is the multiplier ideal sheaf associated to $h$. For a proof we refer to [31] or [8] and [9]. If the metric is smooth then $B$ is a positive line bundle and $\mathscr{I}(h) = \mathcal{O}_X$, the Nadel vanishing theorem reduces to the Kodaira vanishing theorem. A natural question is whether the Nadel vanishing theorem has a generalization for arbitrary Dolbeault cohomology groups. Since a general pseudo-effective line bundle is not necessary semipositive, the critical trick mentioned above could neither be used, and we could not combine it with the approximation technique used by Demailly in [7] to produce such a generalization. Therefore it is difficult to get a vanishing theorem for multiplier ideal sheaf for arbitrary Dolbeault cohomology if the line bundle is only pseudo-effective.

# 3 $(k,s)$-Positivity and Vanishing Theorems

The relation between the Nakano positivity and the Griffiths positivity was studied by Demailly and Skoda in [10], see also [7, 9]. In [7] Demailly introduced the notion of $s$-positivity for any integer $1 \leq s \leq r$ for a vector bundle $E$ of rank $r$. We cited it as the *Demailly s-positivity*. In particular, the Demailly 1-positivity is just the Griffiths positivity and the Demailly $s$-positivity for $s \geq \min\{r, n\}$ is exactly the Nakano positivity. Our aim in this section is to introduce a new concept of positivity for holomorphic vector bundle of arbitrary rank, which unifies the concepts of $k$-positivity for line bundles and the Demailly $s$-positivity for higher rank vector bundles. Then we will give some new vanishing theorems.

Recall the following concept which has been used by Demailly.

**Definition 3.1.** ([7]) A tensor $u \in TX \otimes E$ is called *rank $s$* if $s$ is the smallest non-negative integer such that $u$ can be written as

$$u = \sum_{j=1}^{s} \xi^j \otimes v^j, \quad \xi^j \in TX, \quad v^j \in E.$$

$E$ is said to be *Demailly s-positive* if $i\Theta_h(E)(u,u) > 0$ for any nonzero $u \in TX \otimes E$ of rank $\leq s$. By definition, Demailly 1-positive is equivalent to Griffiths positive and Demailly $s$-positive for $s \geq \{r, n\}$ is equivalent to Nakano positivity. In Demailly's definition, the set of tensors of rank no more than $s$ in $E \otimes TX$ don't form a linear subspace since they are not closed under addition, which causes some difficulties in understanding this notion. To capture the essence of it, we adopt following formulation, which uses the theory of the hermitian form and is a more general definition.

**Definition 3.2.** A holomorphic vector bundle $E$ of rank $r$ with hermitian metric $h$ on a compact complex manifold $X$ of complex dimension $n$ is called $(k,s)$-*positive* for $1 \leq s \leq r$, if the following holds for any $x \in X$ : For any $s$-tuple vectors $v^j \in V, 1 \leq j \leq s$, where $V = E_x$ (resp. $T_xX$), the hermitian form

$$Q_x(\bullet, \bullet) = i\Theta_h(E)(\sum_{j=1}^{s} \cdot \otimes v^j, \sum_{j=1}^{s} \cdot \otimes v^j), \quad \bullet \in W^{\oplus s}, \quad \cdot \in W$$

defined on $W^{\oplus s}$ is semipositive and the dimension of its kernel is at most $k$, where $W = T_xX$ (resp. $E_x$). In this case we write as

$$i\Theta_h(E) >_{(k,s)} 0 \quad \text{or} \quad E >_{(k,s)} 0.$$

Clearly the $(0,s)$-positivity is equivalent to the Demailly $s$-positivity and the Nakano positivity is equivalent to the $(0,s)$-positivity if $s \geq \min\{n,r\}$. The $(0,1)$-positivity is equivalent to the Griffiths positivity. For general integer $k$, the $(k,1)$ positivity is a semipositive version of the Griffiths positivity. A holomorphic vector bundle $E$ of arbitrary rank is called *Griffiths k-positive* if if it is $(k,1)$-positive.

**Example 3.3.** Let $V$ be a complex vector space of dimension $n$ and let $\mathrm{Gr}(V,d)$ denote the Grassmannian of subspaces $W$ of $V$ of codimension $d$. In the following we will prove that the holomorphic tangent bundle of $\mathrm{Gr}(V,d)$ is $((d-1)(n-d-1),1)$-positive. We denote the holomorphic tangent bundle $T\mathrm{Gr}(V,d)$ of $\mathrm{Gr}(V,d)$ by $E$ for brevity. Note that the tangent space at any point is identified with $\mathrm{Hom}(V/W,V)$, which can be seen as the space consisting of $d \times (n-d)$ rectangular complex matrices. We use indices $i, j, k, \cdots$ in the range 1 to $d$; $r, s, t, \cdots$ in the range 1 to $n-d$. Using the homogeneous metric (cf. [5]) and a suitable coordinate system, we have $h_{ir\overline{js}} = \delta_{ij}\delta_{rs}$ and

$$R_{ir\overline{js}kt\overline{lu}} = \delta_{ij}\delta_{kl}\delta_{ru}\delta_{st} + \delta_{il}\delta_{jk}\delta_{rs}\delta_{tu}.$$

For $u = \xi \otimes v = \xi^{ir}\frac{\partial}{\partial z_{ir}} \otimes v^{js}\frac{\partial}{\partial z_{js}} \in T\mathrm{Gr}(V,d) \otimes E$, we have

$$\begin{aligned} i\Theta_h(E)(u,u) &= i\Theta_h(E)(\xi^{ir}\tfrac{\partial}{\partial z_{ir}} \otimes v^{js}\tfrac{\partial}{\partial z_{js}}, \xi^{ir}\tfrac{\partial}{\partial z_{ir}} \otimes v^{js}\tfrac{\partial}{\partial z_{js}}) \\ &= (\delta_{ij}\delta_{kl}\delta_{ru}\delta_{st} + \delta_{il}\delta_{jk}\delta_{rs}\delta_{tu})\xi^{ir}v^{kt}\overline{\xi^{js}}\overline{v^{lu}} \\ &= \sum_{i,k,r,s} \xi^{ir}v^{js}\overline{\xi^{is}}\overline{v^{jr}} + \xi^{ir}v^{js}\overline{\xi^{jr}}\overline{v^{is}} \\ &= \sum_{i,j}|\sum_r \xi^{ir}\overline{v^{jr}}|^2 + \sum_{r,s}|\sum_i \xi^{ir}\overline{v^{is}}|^2. \end{aligned}$$

Fix $(v^{ir})_{d\times(n-d)} \neq 0$. If $i\Theta_h(E)(u,u) = 0$, then the rows and columns of $(\xi^{ir})_{d\times(n-d)}$ are orthogonal to the rows and columns of $(v^{ir})_{d\times(n-d)}$ respectively, hence the rank of the Kernel of $i\Theta_h(E)(\cdot\otimes v,\cdot\otimes v)$ is at most $(d-1)(n-d-1)$. Therefore $E$ is $((d-1)(n-d-1),1)$-positive.

In particular, if $d = 1$, $\mathrm{Gr}(V,1)$ is the $n-1$-dimensional complex projective space $\mathbb{P}^{n-1}$. So the holomorphic tangent bundle of complex projective space is Griffiths positive. Similar calculations as above indicate the holomorphic tangent bundle of complex projective space $\mathbb{P}^{n-1}$ is Nakano semipositive and $((n-1)(s-1),s)$-positive for any $1 \leq s \leq n-1$.

**Proposition 3.4.** *Let $f$ be a surjective holomorphic map from a complex manifold $Y$ of dimension $m$ to a complex manifold $X$ of dimension $n$, and $E$ a $(k,s)$-positive holomorphic vector bundle on $X$ with hermitian metric $h$. Then the pull-back bundle $f^*E$ is $(k+m-n,s)$-positive. In particular, if $f$ is a finite cover then $E$ is $(k,s)$-positive if and only if $f^*E$ is $(k,s)$-positive.*

*Proof.* Let $(e_\alpha)_\alpha$ be a basis of $E$, then $(f^*e_\alpha)_\alpha$ is a basis of $f^*E$. Note $i\Theta_{f^*h}(f^*E) = if^*\Theta_h(E)$ and

$$i\Theta_{f^*h}(f^*E)(\sum_i(\cdot \otimes f^*e_{\alpha_i}), \sum_i(\cdot \otimes f^*e_{\alpha_i})) = \Theta_h(E)(\sum_i(f_*(\cdot) \otimes e_{\alpha_i}), \sum_i(f_*(\cdot) \otimes e_{\alpha_i}))$$

So the hermitian form $i\Theta_{f^*h}(f^*E)(\sum_i(\cdot \otimes f^*e_{\alpha_i}), \sum_i(\cdot \otimes f^*e_{\alpha_i}))$ is semipositive and the rank of its kernel is at most $m - n + k$. Thus $f^*E$ is $(m - (n-k), s)$-positive. $\square$

**Proposition 3.5.** *If $0 \to S \to E \to Q \to 0$ is an exact sequence of hermitian vector bundles on a complex manifold $X$ of dimension $n$, then for any integer $1 \le k \le n$, we have*
a). $E <_{(k,s)} 0 \Longrightarrow S <_{(k,s)} 0$ *for* $s \le \text{rank } S$;
b). $E >_{(k,1)} 0 \Longrightarrow Q >_{(k,1)} 0$;
c). *If the exact sequence splits, then* $S >_{(k,s)} 0$ *and* $Q >_{(k,s)} \Longleftrightarrow E >_{(k,s)} 0$ *for* $s \le \min\{\text{rank } S, \text{rank } Q\}$.

*Proof.* Let $\{e_\alpha\}$ and $\{e_\mu\}$ be the basis of $S$ and $Q$ respectively, and let $\{e^\beta\}$ and $\{e^\nu\}$ be the dual basis of $\{e_\alpha\}$ and $\{e_\mu\}$ respectively. Let
$$B = B^\mu_{\alpha j} e^\alpha \otimes e_\mu \otimes dz^j \in \Omega^{1,0}(\text{Hom}(S,Q))$$
be the second fundamental form associated with the exact sequence. Then
$$B^* = \sum_{\mu,\alpha,j} \overline{B^\mu_{\alpha j}} e^\mu \otimes e_\alpha \otimes d\bar{z}^j \in \Omega^{0,1}(\text{Hom}(Q,S)).$$

With the induced hermitian metric of $E$, the Chern curvature form of $S$ and $Q$ are (cf. [17, 7]):
$$\begin{aligned} i\Theta(S) &= i\Theta(E)|_S + B^* \wedge B; \\ i\Theta(Q) &= i\Theta(E)|_Q + B \wedge B^*. \end{aligned} \qquad (3.1)$$

Thus for $u = \sum_{\alpha j} u^{\alpha j} \frac{\partial}{\partial z_j} \otimes e_\alpha \in TX \otimes S$ and $v = \sum_{\mu k} v^{\mu k} \frac{\partial}{\partial z_k} \otimes e_\mu \in TX \otimes Q$, we have
$$\begin{aligned} i\Theta(S)(u,u) &= i\Theta(E)|_S(u,u) - \sum_{\mu,\nu,\alpha,\beta,j,k} h_{\mu\bar\nu} B^\mu_{\alpha j} \overline{B^\nu_{\beta k}} u^{\alpha j} \overline{u^{\beta k}}; \\ i\Theta(Q)(v,v) &= i\Theta(E)|_Q(v,v) + \sum_{\mu,\nu,\alpha,\beta,j,k} h_{\alpha\bar\beta} \overline{B^\mu_{\alpha j}} B^\nu_{\beta k} v^{\mu k} \overline{v^{\nu j}}. \end{aligned}$$

Thus
$$i\Theta(S)(u,u) = i\Theta(E)|_S(u,u) - |\sum_{\alpha,j} B^\bullet_{\alpha j} u^{\alpha j}|^2, \qquad (3.2a)$$
where $B^\bullet_{\alpha j} u^{\alpha j} = \sum_\mu B^\mu_{\alpha j} u^{\alpha j} e_\mu \in Q$. From (3.2a) we conclude a). If $v = (\sum_k \xi^k \frac{\partial}{\partial z_k}) \otimes (\sum_\mu w^\mu e_\mu)$ is a tensor of rank 1, then
$$i\Theta(Q)(u,u) = i\Theta(E)|_Q(u,u) + |\sum_{\mu,j} \overline{B^\mu_{\bullet j}} w^\mu \overline{\xi^j}|^2, \qquad (3.2b)$$
where $\overline{B^\mu_{\bullet j}} w^\mu \overline{\xi^j} = \sum_\alpha \overline{B^\mu_{\alpha j}} w^\mu \overline{\xi^j} e_\alpha$. By (3.2b) we have b). Note that c) is clearly true since $B = 0$ in this case. $\square$

**Proposition 3.6.** *Let $E$ and $F$ be hermitian holomorphic vector bundles.*
a). $E >_{(k,s)} 0 \Longrightarrow E >_{(k+1,s)} 0$, *and* $E >_{(k,s-1)} 0$ *if* $s \ge 2$;
b). $E >_{(a,p)} 0$ *and* $F >_{(b,q)} 0 \Longrightarrow E \otimes F >_{(c,s)} 0$ *with* $c = \max\{a,b\}$ *and* $s = \min\{p,q\}$;
c). $E >_{(k,s)} 0 \Longrightarrow \Gamma^a E >_{(k,s)} 0$, *where $\Gamma^a E$ is an irreducible representation of $Gl(E)$ of highest weight $a = (a_1, a_2, \cdots, a_r) \in \mathbb{Z}^r$, with $a_1 \ge a_2 \ge \cdots \ge a_r \ge 0$. In particular, for $p$-order symmetric product $S^p E$ and $q$-order wedge $\wedge^q E, (1 \le q \le \text{rank } E)$,*

*Proof.* a) follows by definitions. Let $h_E, h_F$ be the corresponding hermitian metric on $E$ and $F$ with Chern curvature forms $R^E, R^F$ respectively, and $(t_j)_j, (e_\alpha)_\alpha, (f_\beta)_\beta$ be the orthonormal frames of $TX, E, F$ respectively, let $u = u^{\alpha\mu j} t_j \otimes e_\alpha \otimes f_\mu \in E \otimes F \otimes TX$ be a tensor of rank at most $s$. Since $\Theta_{h_E \otimes h_F}(E \otimes F) = \Theta_{h_E}(E) \otimes Id_F + Id_E \otimes \Theta_{h_F}(F)$, we have

$$\begin{aligned}
& i\Theta_{h_E \otimes h_F}(E \otimes F)(u,u) \\
=\ & \sum_{\alpha,\beta,j,k,\mu} R^E_{\alpha\bar{\beta}j\bar{k}} u^{\alpha\mu j}\overline{u^{\beta\mu k}} + \sum_{\mu,\nu,j,k,\alpha} R^E_{\mu\bar{\nu}j\bar{k}} u^{\alpha\mu j}\overline{u^{\alpha\nu k}} \\
=\ & \sum_\mu (\sum_{\alpha,\beta,j,k} R^E_{\alpha\bar{\beta}j\bar{k}} u^{\alpha\mu j}\overline{u^{\beta\mu k}}) + \sum_\alpha (\sum_{\mu,\nu,j,k} R^E_{\mu\bar{\nu}j\bar{k}} u^{\alpha\mu j}\overline{u^{\alpha\nu k}}).
\end{aligned}$$

Fix $\mu$ and $\alpha$ respectively, consider $u^{\alpha\mu j} t_j \otimes e_\alpha \otimes f_\mu$ as the tensors of $E \otimes TX$ and $F \otimes TX$ respectively, of rank at most $s = \min\{p,q\}$, we conclude $E \otimes F >_{(c,s)} 0$ for $c = \max\{a,b\}$. This proves b). If $E >_{(k,s)} 0$, using b) inductively we have $E^{\otimes p} >_{(k,s)} 0$. Since $S^p E$ and $\wedge^q E$ are direct summand of $E^{\otimes p}$ and $E^{\otimes q}$ respectively, by c) of Proposition 3.5 we have $S^p E >_{(k,s)} 0$ and $\wedge^q E >_{(k,s)} 0$. Using b) again we have

$$S^{a_1} E \otimes S^{a_2} E \otimes \cdots \otimes S^{a_r} E >_{(k,s)} 0.$$

Since $\Gamma^a E$ is a direct summand in $S^{a_1} E \otimes S^{a_2} E \otimes \cdots \otimes S^{a_r} E$, we have $\Gamma^a E >_{(k,s)} 0$. $\square$

**Proposition 3.7.** *If $E$ is a hermitian holomorphic line bundle, then $E$ is $k$-positive if and only if $E$ is Griffiths $k$-positive.*

*Proof.* If $r = \text{rank } E = 1$, then the rank 1-tensors of $E \otimes TX$ is identified with $TX$, thus $E$ is Griffiths $k$-positive if and only if $i\Theta_h(E)$ is semipositive on $TX$ and has at least $n - k$-positive eigenvalues, *i.e.*, if and only if $E$ is $k$-positive. $\square$

For a holomorphic vector bundle $E$ of higher rank $r > 1$, we define the projective bundle

$$\pi : \mathbb{P}(E) \to X,$$

whose fibre $\pi^{-1}(x)$ at $x \in X$ is $(r-1)$-dimensional subspace of $E|_x$. Then $\mathbb{P}(E)$ carries a tautological subbundle $F$ of $\pi^* E$ given by

$$F|_y = y \subset (\pi^* E)_y = E|_{\pi(y)}, \quad for \quad y \in \mathbb{P}(E).$$

We define a holomorphic vector bundle $\xi_E$ on $\mathbb{P}(E)$ by

$$\xi_E = \pi^* E / F.$$

Recall that a vector bundle $E$ is called *$k$-positivity* if $\xi_E$ is $k$-positive (cf.[40, 20]).

**Theorem 3.8.** *If $E$ is a Griffiths $k$-positive holomorphic vector bundle then it is $k$-positive. In particular if $E$ is a Griffiths $k$-positive holomorphic vector bundle on a compact Kähler manifold $X$ of complex dimension $n$. Then $H^{p,q}(X,E) = 0$ for $p + q \geq n + k + r$.*

*Proof.* Let $h$ be the corresponding hermitian metric on $E$. Note that the dual bundle $\xi_E^* = (\pi^* E/F)^*$ is isomorphic to the conormal bundle $F^\perp \subset \pi^* E^*$ of $F$. For any $v \in \xi_E^* - \{0\}$, there is a naturally associated hermitian metric on $\xi_E^*$ defined by $H(z,v)_\omega = h_{j\bar{k}} v^j \bar{v}^k = (v,v)_z$ for $(z,[v]) \in \mathbb{P}(E^*)$. Then $H(z, \lambda v) = |\lambda|^2 H(z,v)$ for any $\lambda \in \mathbb{C}^*$. Hence $H$ defines a hermitian metric on $\xi_E^*$. Choose normal coordinates $z^1, \cdots, z^n$ such that $h(0) = I$ and $dh(0) = 0$. Then $\Theta_H(\xi_E^*) = -\partial\bar{\partial} \log H(z,v) = -\partial\bar{\partial} \log(h_{j\bar{k}} v^j \bar{v}^k)$, and

$$\Theta_H(\xi_E^*) = \frac{\Theta_h(E^*)(v,v)}{|v|^2} - \frac{(\partial v, \partial v)|v|^2 - (\partial v, v)(v, \partial v)}{|v|^4}.$$

Therefore
$$\begin{aligned}\Theta_H(\xi_E) &= \frac{\Theta_h(E)(v,v)}{|v|^2} + \frac{(\partial v, \partial v)|v|^2 - (\partial v, v)(v, \partial v)}{|v|^4} \\ &= -\frac{(\partial\bar{\partial} h_{j\bar{k}}) v^j \bar{v}^k}{h_{j\bar{k}} v^j \bar{v}^k} + \frac{h_{j\bar{k}} \partial v^j \wedge \overline{\partial v}^k}{h_{j\bar{k}} v_j \bar{v}_k} - \frac{h_{j\bar{k}} \partial v^j \bar{v}^k \wedge h_{j\bar{k}} v^j \overline{\partial v}^k}{(h_{j\bar{k}} v^j \bar{v}^k)^2}.\end{aligned}$$

In matrix form, the curvature

$$\Theta_H(\xi_E)|_{(z,v)} = \begin{pmatrix} \left(\frac{\Theta_h(E)(v,v)}{|v|^2}\right)_{n\times n} & \bigg| & 0 \\ \rule{2cm}{0.4pt} & & \rule{3cm}{0.4pt} \\ 0 & \bigg| & \left(\frac{(1+|w|^2)\delta_{jk} - w_k \bar{w}_j}{(1+|w|^2)^2}\right)_{(r-1)\times(r-1)} \end{pmatrix},$$

where $w = (\frac{v^0}{v^j}, \cdots, \frac{v^{j-1}}{v^j}, \frac{v^{j+1}}{v^j}, \cdots, \frac{v^r}{v^j})$ is the local coordinate of the open subset $V_j = \{v^j \neq 0\}$ of the fiber $\mathbb{P}(E|_z)$. Note that $\left(\frac{(1+|w|^2)\delta_{ij} - w_j \bar{w}_i}{(1+|w|^2)^2}\right)_{(r-1)\times(r-1)}$ is a positive definite matrix with eigenvalues $1/(1+|w|^2)$ of order $r-2$ and eigenvalue $1/(1+|w|^2)^2$ of order 1. Therefore $i\Theta_h(E)$ is semipositive if and only if $i\Theta_H(\xi_E)$ is semipositive. If $E$ is Griffiths $k$-positive, then the rank 1-tensors of $E \otimes TX$ is identified with $TX \cup E$, which in particular implies that $\frac{\Theta_h(E)(v,v)}{|v|^2}$ is semipositive and has at least $n-k$-positive eigenvalues for any nonzero $v \in E$. Then $i\Theta_H(\xi_E)$ has at least $n-k+r-1$-positive eigenvalues. Thus if $E$ is Griffiths $k$-positive then $\xi_E$ is $k$-positive.

Let $E$ be a Griffiths $k$-positive holomorphic vector bundle on the compact Kähler manifold $X$. By the Gigante-Girbau vanishing theorem, when $p + q > (n + r - 1) + k$, then $H^p(\mathbb{P}(E), \Omega^q_{\mathbb{P}(E)}(\xi_E)) = 0$. Using the Griffiths-Le Potier-Schneider isomorphism (See Remark 5.9 of Section 5):

$$H^q(\mathbb{P}(E), \Omega^p_{\mathbb{P}(E)}(\xi_E)) \cong H^q(X, \Omega^p_X(E)), \quad n \geq p, q \geq 0, \tag{3.3}$$

we have $H^p(X, \Omega^q_X(E)) = 0$ for $p + q \geq n + r + k$. $\square$

Recall from (2.10) we have

$$\langle [i\Theta_h(E), \Lambda]u, u \rangle_\omega = \sum_{S_{q-1}, \alpha, \beta, j, k} R_{\alpha\bar{\beta}j\bar{k}} u^\alpha_{\bar{j}\bar{S}_{q-1}} \overline{u^\beta_{\bar{k}S_{q-1}}} \tag{3.4}$$

for any $(n,q)$-form $u = \sum u_{\bar{K}}^\alpha dz^1 \wedge \cdots \wedge dz^n \wedge d\bar{z}^K \otimes e_\alpha$ with value in $E$. Since $u_{\bar{j}\bar{S}_{q-1}}^\alpha = 0$ for $\bar{j} \in \bar{S}_{q-1}$, so the rank (in the sense of Demailly) of $(u_{\bar{j}\bar{S}_{q-1}}^\alpha)_{\bar{j},\alpha}$ is in fact no more than $\min\{n-q+1, r\}$.

**Theorem 3.9.** *Let $X$ be a compact Kähler manifold of dimension $n$ and $E$ a hermitian holomorphic vector bundle of rank $r$ on $X$ such that $E >_{(k,s)} 0$. Then*

$$H^{n,q}(X, E) = 0, \quad for \quad q > k \quad and \quad s \geq \min\{n-q+1, r\}.$$

*Proof.* We may assume $k \leq n$. For any $E$-valued $(n,q)$-form $u = \sum u_{\bar{K}}^\alpha dz^1 \wedge \cdots \wedge dz^n \wedge d\bar{z}^K \otimes e_\alpha \in \Omega^{n,q}(E)$, by Definition 3.2, we could diagonalize the hermitian form $i\Theta_h(E)(u, u)$ since it is semipositive. Therefore we could write $\sum_{S_{q-1}, \alpha, \beta, j, k} R_{\alpha\bar{\beta}j\bar{k}} u_{\bar{j}\bar{S}_{q-1}}^\alpha \overline{u_{\bar{k}S_{q-1}}^\beta}$ as

$$\sum_{S_{q-1},\alpha,\beta,j,k} R_{\alpha\bar{\beta}j\bar{k}} u_{\bar{j}\bar{S}_{q-1}}^\alpha \overline{u_{\bar{k}S_{q-1}}^\beta} = \sum_{S_{q-1},\alpha,j} \lambda_j^\alpha u_{\bar{j}\bar{S}_{q-1}}^\alpha \overline{u_{\bar{j}S_{q-1}}^\alpha},$$

where the eigenvalues $(\lambda_j^\alpha)_{1 \leq j \leq n, 1 \leq \alpha \leq r}$ are non-negative. For fixed $\alpha$, without loss of generality, assume that $\lambda_1^\alpha \leq \lambda_2^\alpha \leq \cdots \leq \lambda_n^\alpha$ with $\lambda_{k+1}^\alpha > 0$. Put $\lambda = \min\{\lambda_{k+1}^\alpha | 1 \leq \alpha \leq r\}$. Then $\lambda$ is a positive number. If $q > k$ then

$$\begin{aligned}
\sum_{S_{q-1},\alpha,\beta,j} R_{\alpha\bar{\beta}j\bar{k}} u_{\bar{j}\bar{S}_{q-1}}^\alpha \overline{u_{\bar{k}S_{q-1}}^\beta} &= \sum_{j,\bar{S}_{q-1}} \sum_\alpha \lambda_j^\alpha |u_{\bar{j}\bar{S}_{q-1}}^\alpha|^2 \\
&= \sum_K \sum_\alpha \sum_{j \in K} \lambda_j^\alpha |u_{\bar{K}}^\alpha|^2 \\
&\geq \sum_K \sum_\alpha (\lambda_1^\alpha + \lambda_2^\alpha + \cdots + \lambda_q^\alpha) |u_{\bar{K}}^\alpha|^2 \\
&\geq \lambda (\sum_K \sum_\alpha |u_{\bar{K}}^\alpha|^2) \\
&= \lambda |u|_\omega^2.
\end{aligned}$$

Thus when $q > k$, then $[i\Theta_h(E), \Lambda]$ is positive definite on $E$-valued $(n,q)$-forms. So we have $H^{n,q}(X, E) = 0$ for $q > k$ and $s \geq \min\{n-q+1, r\}$. □

The special case when $k = 0$ of Theorem 3.9 was proved by Demailly in [9].

**Proposition 3.10.** *Let $E$ be a rank-$r$ holomorphic vector bundle with hermitian metric $h$ on a complex manifold $X$ of dimension $n$. If $E$ is Griffiths $k$-positive, then for any integer $1 \leq s \leq \min\{r, n\}$,*

$$i\Theta_h(E) + \text{Tr}_E(i\Theta_h(E)) \otimes h >_{(k,s)} 0. \tag{3.5}$$

*Proof.* By Definition 3.2, it suffices to prove (3.5) for any subbundle $G \otimes S$ of $E \otimes TX$ with $G = E$ and rank $S = s$, or $S = TX$ and rank $G = s$. Let $(t_j)_{1 \leq j \leq m}$ be a basis of $S$ and $(e_\alpha)_{1 \leq \alpha \leq \rho}$ an orthonormal basis of $G$. The coefficients $R_{\alpha\bar{\beta}j\bar{k}}$ of $i\Theta_h(E)$ with respect to the basis $t_j \otimes e_\alpha$ satisfy the symmetric relation $\overline{R_{\alpha\bar{\beta}j\bar{k}}} = R_{\beta\bar{\alpha}k\bar{j}}$. For any $\xi = \xi^j t_j \in S$ and $u = u^{\alpha j} e_\alpha \otimes t_j \in E \otimes TX$, we have the formulae

$$i\Theta_h(E)(u,u) = \sum_{\alpha,\beta,j,k} R_{\alpha\bar{\beta}j\bar{k}} u^{\alpha j} \overline{u^{\beta k}},$$

$$\mathrm{Tr}_E(i\Theta_h(E))(\xi,\xi) = \sum_{\alpha,j,k} R_{\alpha\bar{\alpha}j\bar{k}}\xi^j\overline{\xi^k},$$

$$(i\Theta_h(E) + \mathrm{Tr}_E(i\Theta_h(E))\otimes h)(u,u) = \sum_{\alpha,\beta,j,k} R_{\alpha\bar{\beta}j\bar{k}}u^{\alpha j}\overline{u^{\beta k}} + R_{\alpha\bar{\alpha}j\bar{k}}u^{\alpha j}\overline{u^{\alpha k}}. \tag{3.6}$$

We use the notions and techniques of Lemma 10.15 and Proposition 14 in [9]. Let $\sigma$ describe the set $U_q^\rho$ of $\rho$-tuples of $q$-th roots of unity and put

$$u'^j_\sigma = \sum_{1\leq\alpha\leq\rho} u^{\alpha j}\bar{\sigma}_\alpha \in \mathbb{C},$$

$$\hat{u}_\sigma = \sum_j u'^j_\sigma t_j \in S, \qquad \hat{e}_\sigma = \sum_\alpha \sigma_\alpha e_\alpha \in G.$$

By Lemma 10.15 of [9], we have

$$\begin{aligned}
& q^{-\rho}\sum_{\sigma\in U_q^\rho} i\Theta_h(E)(\hat{u}_\sigma\otimes\hat{e}_\sigma, \hat{u}_\sigma\otimes\hat{e}_\sigma) \\
= & q^{-\rho}\sum_{\sigma\in U_q^\rho} R_{\alpha\bar{\alpha}j\bar{k}}u'^j_\sigma \bar{u}'^k_\sigma \sigma_\alpha\bar{\sigma}_\beta \\
= & \sum_{j,k,\alpha\neq\beta} R_{\alpha\bar{\beta}j\bar{k}}u^{\alpha j}\overline{u^{\beta k}} + R_{\alpha\bar{\alpha}j\bar{k}}u^{\alpha j}\overline{u^{\alpha k}}.
\end{aligned} \tag{3.7}$$

Clearly we have

$$\sum_{\alpha,k,j} R_{\alpha\bar{\alpha}j\bar{k}}u^{\alpha j}\overline{u^{\alpha k}} = \sum_\alpha i\Theta_h(E)((\sum_j u^{\alpha j}t_j)\otimes e_\alpha, (\sum_j u^{\alpha j}t_j)\otimes e_\alpha). \tag{3.8}$$

By (3)- (3.8) we get

$$\begin{aligned}
& (i\Theta_h(E) + \mathrm{Tr}_E(i\Theta_h(E))\otimes h)(u,u) \\
= & q^{-\rho}\sum_{\sigma\in U_q^\rho} i\Theta_h(E)(\hat{u}_\sigma\otimes\hat{e}_\sigma, \hat{u}_\sigma\otimes\hat{e}_\sigma) \\
& + \sum_\alpha i\Theta_h(E)((\sum_j u^{\alpha j}t_j)\otimes e_\alpha, (\sum_j u^{\alpha j}t_j)\otimes e_\alpha).
\end{aligned}$$

Since $i\Theta_h(E) >_{(k,1)} 0$, if we fix the $E$-components of the tensors in the summation of the right hand side, we get a summation of semipositive definite matrices of rank at least $n-k$, and hence the left hand side is also semipositive definite of rank at least $n-k$ if we fix the $E$-components, thus $i\Theta_h(E) + \mathrm{Tr}_E(i\Theta_h(E))\otimes h >_{(k,s)} 0$. □

**Proposition 3.11.** *Let $E$ be a holomorphic vector bundle of rank $r \geq 2$ with hermitian metric $h$ on a complex manifold $X$ of dimension $n$. If $E$ is Griffiths $k$-positive, then for any integer $1\leq s\leq \max\{r,n\}$, we have*

$$s\mathrm{Tr}_E(i\Theta_h(E))\otimes h - i\Theta_h(E) >_{(k,s)} 0. \tag{3.9}$$

*Proof.* Let us distinguish two cases.

*Case* a). $m = 1$. Let $u \in E \otimes TX$ be a tensor of rank 1. Then $u$ can be written as $u = \xi_1 \otimes e_1$ with $\xi_1 \in TX, \xi_1 \neq 0$, and $e_1 \in E, |e_1| = 1$. Complete $e_1$ into an orthonormal basis $(e_1, \cdots, e_r)$ of $E$. One gets immediately

$$\begin{aligned}(\mathrm{Tr}_E(i\Theta_h(E)) \otimes h)(u,u) &= \sum_\alpha i\Theta_h(E)(\xi_1 \otimes e_\alpha, \xi_1 \otimes e_\alpha) \\ &\geq i\Theta_h(E)(\xi_1 \otimes e_1, \xi_1 \otimes e_1) = i\Theta_h(E)(u,u).\end{aligned} \quad (3.10)$$

*Case* b). $m \geq 2$. As in the proof of Proposition 3.7, it suffices to prove (3.9) for subbundles $G \otimes TX$ of $E \otimes TX$ with rank $G = s$. Let $(e_\alpha)_{1 \leq \alpha \leq s}$ be an orthonormal basis of $G$ which complete to a basis $(e_\alpha)_{1 \leq \alpha \leq r}$ of $E$, and let $i\Theta_h(G)$ be the restriction of $i\Theta_h(E)$ to $G \otimes TX$. Then *Case* a) shows that

$$\mathrm{Tr}_G(i\Theta_h(G)) \otimes h - i\Theta_h(G) >_{(k,1)} 0.$$

Note that in the Proposition 3.10, we could substitute $i\Theta_h(E)$ and $h$ by any hermitian form on $E \otimes TX$ and $E$ respectively, then the conclusion of Proposition 3.10 is still right in the general case. Consider $\Theta' := \mathrm{Tr}_G(i\Theta_h(G)) \otimes h - i\Theta_h(G)$ as an hermitian form on $TX \otimes G$ and restrict $h$ to $G$ as a hermitian form, using Proposition 3.10 we get

$$\Theta' + \mathrm{Tr}_G \Theta' \otimes h = s\mathrm{Tr}_G(G) \otimes h - \Theta_h(G) >_{(k,s)} 0.$$

Vary $G$ as a rank $s$-subbundle of $E$ then we have $s\mathrm{Tr}_E(E) \otimes h - \Theta_h(E) >_{(k,s)} 0$. □

**Corollary 3.12.** *Let $E$ be a Griffiths $k$-positive hermitian holomorphic vector bundle of rank $r \geq 2$. Then for any integer $m \geq 1$,*

$$E^* \otimes (\det E)^s >_{(k,s)} 0.$$

*Proof.* Apply Proposition 3.11 to $E^* \otimes TX$ and note that

$$i\Theta(E^* \otimes (\det E)^s) = s(i\Theta(\det E)) \otimes h + i\Theta(E^*) = s\mathrm{Tr}_E(i\Theta(E)) \otimes h - i\Theta(E)^t$$

and $\mathrm{Tr}_E(i\Theta(E)) = \mathrm{Tr}_E(i\Theta(E)^t)$. □

Theorem 3.9 in combination with Proposition 3.10 and Corollary 3.12 immediately imply the following consequences which generalize (10.19) in [9]:

**Corollary 3.13.** *Let $E$ be a Griffiths $k$-positive hermitian holomorphic vector bundle of rank $r \geq 2$ on a compact Kähler manifold of dimension $n$. Then for any integer $s \geq 1$, we have*
a). $H^q(X, K_X \otimes E \otimes \det E) = 0 \quad \text{for} \quad q > k$;
b). $H^q(X, K_X \otimes E^* \otimes (\det E)^s) = 0 \quad \text{for} \quad q > k \quad \text{and} \quad s \geq \min\{n-q+1, r\}$.

The following generalized Lefschetz theorem (the classical Lefschetz theorem is the special case when $r=1$ and $k=0$) is a generalization of theorem H in [17], where $E$ was supposed to be Griffiths positive.

**Theorem 3.14.** *Let $(E, h)$ be a Griffiths $k$-positive vector bundle on a compact complex manifold $X$ and $v \in H^0(X, E)$ a holomorphic section whose zero locus $S \subset X$ is a smooth submanifold of codimension $r$. Then we have $H_q(S, \mathbb{Z}) \to H_q(X, \mathbb{Z}) \to 0$ for $q \leq n - r - k$ and $0 \to H_q(S, \mathbb{Z}) \to H_q(X, \mathbb{Z}) \to 0$ for $q \leq n - r - k - 1$.*

*Proof.* Let $v(z) = (v^1(z), \cdots, v^r(z))$. Define a function
$$f(z) = h_{\alpha\bar\beta}(z) v^\alpha(z) \overline{v^\beta(z)}.$$
Choose local coordinates $z^1, \cdots, z^n$ on $X$ such that $v(z) = (z^1, \cdots, z^r)$. We may assume $h(0) = I$. Then $f(z) = \sum_{\alpha=1}^r |z^\alpha|^2$, and $df(0) = 0$ and the Hessian of $f$ at origin is
$$\mathrm{Hess}(f) = \begin{pmatrix} I_r & 0 \\ 0 & 0 \end{pmatrix}.$$
So $S$ is a non-degenerate critical submanifold of $f$.

Let $z_0$ be a critical point of $f$ on $X - S$. Recall that the index of $f$ at $z_0$ is the dimension of the subspace of $T_{z_0} X$ on which $\mathrm{Hess}(f)$ is negative definite. We claim that the index of $f$ is no less than $n - r - k + 1$.

In fact, we may choose normal coordinate such that $h(z_0) = I$ and $dh(z_0) = 0$. Since $z_0$ is a critical point of $f$, we have
$$\sum_\alpha (\overline{v^\alpha} dv^\alpha + v^\alpha \overline{dv^\alpha}) = 0, \tag{3.11}$$
and $\partial\bar\partial f = ((\partial\bar\partial h_{\alpha\bar\beta}) v^\alpha \overline{v^\beta} + h_{\alpha\bar\beta} dv^\alpha \wedge \overline{dv^\beta}$. Thus
$$(\mathrm{Hess}(f)(z_0))_{j\bar k} = -\sum_{\alpha\beta} R^\beta_{\alpha j \bar k} v^\alpha \overline{v^\beta} + \sum_\alpha \frac{\partial v^\alpha}{\partial z^j} \overline{\frac{\partial v^\alpha}{\partial z^k}}.$$
Thus for $t = t^j \frac{\partial}{\partial z_j} \in T_{z_0} X$, we have
$$\mathrm{Hess}(f)_{j\bar k} t^j \overline{t^k} = -i\Theta_h(E)(t \otimes v, t \otimes v) + \sum_\alpha |\sum_j \frac{\partial v^\alpha}{\partial z^j} t^j|^2.$$
Let
$$W = \{t \in T_{z_0} X \mid \sum_j \frac{\partial v^\alpha}{\partial z^j} t^j = 0, 1 \leq \alpha \leq r\}.$$
Then $\dim W = n - \mathrm{rank}\,(\frac{\partial v^\alpha}{\partial z^j})_{r \times n}$. Note by (3.11) we have $\sum_j \frac{\partial v^\alpha}{\partial z^j} \overline{v^\alpha} = 0$. Thus $\mathrm{rank}\,(\frac{\partial v^\alpha}{\partial z^j})_{r \times n} \leq r - 1$. So $\dim W \geq n - r + 1$. Note the dimension of the kernel of $-i\Theta_h(E)(\cdot \otimes v, \cdot \otimes v)$ is at most $k$ since $E$ is Griffiths $k$-positive. Thus the index of $f$ is at least $n - r + 1 - k$. The same argument as the proof of theorem H in [17] by using the Morse theory implies the conclusions of Theorem Theorem 3.14. $\square$

Let $\mathscr{I}_S$ be the ideal sheaf of $S$. If $X$ is a compact Kähler manifold, then $S$ is a Kähler submanifold. In this case we have the following commutative diagram

$$\begin{array}{ccccccc}
& & 0 & & 0 & & \\
& & \downarrow & & \downarrow & & \\
\cdots \longrightarrow H^q(\mathscr{I}_S) \longrightarrow & & H^q(\mathscr{O}_X) \longrightarrow & & H^q(\mathscr{O}_S) \longrightarrow & & H^{q+1}(\mathscr{I}_S) \\
& & \downarrow & & \downarrow & & \\
& & H^q(X,\mathbb{C}) \longrightarrow & & H^q(S,\mathbb{C}). & &
\end{array}$$

Using that $0 \to H^q(X,\mathbb{C}) \to H^q(S,\mathbb{C})$ for $q \leq n-r-k$ and $0 \to H^q(X,\mathbb{C}) \to H^q(S,\mathbb{C}) \to 0$ for $q \leq n-r-k-1$, we get the following vanishing cohomology

$$H^q(\mathscr{I}_S) = 0, \quad for \ \ 0 \leq q \leq n-r-k. \tag{3.12}$$

If $E$ is a line bundle $B$, then $S$ is a hypersurface, we get $H^q(X, B^*) = 0$ for $0 \leq q \leq n-k-1$. We have a pair of short exact sequences

$$0 \longrightarrow \Omega_X^p(B^*) \longrightarrow \Omega_X^p \longrightarrow \Omega_X^p|_S \longrightarrow 0$$

$$0 \longrightarrow \Omega_S^{p-1}(B^*) \longrightarrow \Omega_X^p|_S \longrightarrow \Omega_S^p \longrightarrow 0 \ .$$

In cohomology, we have:

$$\begin{array}{ccccccc}
H^q(X, \Omega_X^p(B^*)) & \longrightarrow & H^q(X, \Omega_X^p) & \longrightarrow & H^q(S, \Omega_X^p|_S) & \longrightarrow & \cdots \\
 & & & & \| & & \\
H^{q-1}(S, \Omega_S^{p-1}) & \longrightarrow & H^q(S, \Omega_S^{p-1}(B^*)) & \longrightarrow & H^q(S, \Omega_X^p|_S) & \longrightarrow & H^q(S, \Omega_S^{p-1}).
\end{array}$$

Using that the map $H^q(X,\mathbb{C}) \to H^q(S,\mathbb{C})$ is injective for $q \leq n-r-k$ and bijective for $q \leq n-r-k-1$, we get

$$H^q(X, \Omega_X^p(B^*)) = 0, \quad for \ \ 0 \leq p+q \leq n-1-k,$$

which is equivalent to the Gigante-Girbau vanishing theorem by the Serre duality theorem.

## 4 Sommese's Vanishing theorems

Sommese introduced the concept of $k$-ample holomorphic vector bundles on a compact complex manifold in [42]. A holomorphic line bundle $B$ on a compact complex manifold $X$ is called *k-ample* if some powers $B^m$ is spanned and the maximum of the dimensions of the fibres of the natural map $\Phi_{B^m}$ associated with $H^0(X, B^m)$ is at most $k$. Suppose that $s_0, s_1, \cdots, s_N$ is a basis of $H^0(X, B^m)$, the holomorphic map $\Phi_{B^m} : X \to \mathbb{P}^N$ is defined by sending $x \in X$ to $[s_0(x) : s_1(x) : \cdots : s_N(x)]$. A holomorphic vector bundle $E$ of arbitrary rank is called $k$-*ample* if the tautological line bundle $\xi_E$ on $\mathbb{P}(E)$ is $k$-ample. Sommese established some vanishing theorems when $X$ is projective using algebraic geometric tools. An interesting phenomena is that the vanishing theorems obtained for $k$-ample vector bundles hold for $k$-positive vector bundles on compact Kähler manifolds as well. Naturally, we conjecture that $k$-ampleness is equivalent to $k$-positivity, at least on compact Kähler manifold. However, even on a projective algebraic manifold a $k$-positive line bundle is not necessary $k$-ample, as indicated by the following example:

**Example 4.1.** Let $\mathscr{E}$ be an elliptic curve and $A, B, C \in \mathscr{E}$ three different points. Let $B_1 = \mathscr{O}(A - B)$ and $B_2 = \mathscr{O}(C)$ be two different holomorphic line bundles. Then $B_2$ is ample hence a positive line bundle. Let $h_1$ be a hermitian metric on $B_1$ such that $i\Theta_{h_1}(B_1) = 0$ and $h_2$ a hermitian metric on $B_2$ such that $i\Theta_{h_2}(B_2) > 0$. Let $X$ be the algebraic surface $\mathscr{E} \times \mathscr{E}$ and $B$ be the line bundle $\pi_1^* B_1 \otimes \pi_2^* B_2$ on $X$ with the hermitian metric $h = \pi_1^* h_1 \times \pi_2^* h_2$. Then $B$ is $(1,1)$-positive and hence 1-positive. But $H^0(X, B^n) = \pi_1^* H^0(\mathscr{E}, \mathscr{O}(nA - nB)) \otimes \pi_2^* H^0(\mathscr{E}, \mathscr{O}(nC)) = \{0\}$. Hence $B$ is not spanned and hence is not $k$-ample for any $k$.

In the book [40], the authors said that there is no explicit relations between the concepts of $k$-ampleness and $k$-positivity. However, in 1983, Sommese proved the following Lemma (cf. [43], Lemma (1.1)) for any holomorphic vector bundle on a compact complex manifold, form which that $k$-ampleness implying $k$-positivity follows easily.

**Lemma 4.2.** *Let $\psi$ be a holomorphic map between a compact complex manifold $X$ and a hermitian manifold $Y$ with hermitian metric $h$. Assume the maximum of the dimensions of the fibres of $\psi$ is at most $k$. Then there exist a smooth non-negative function $f$ such that $\partial\bar{\partial}f + \psi^* h$ is semipositive and has at least $n - k$ positive eigenvalues.*

**Theorem 4.3.** *If $E$ is a $k$-ample line bundle on a compact complex manifold $X$, then it is $k$-positive.*

*Proof.* It suffices to prove it when $E$ is a line bundle, denoted by $B$. From the definition of $k$-ampleness, there exists some positive integer $m$ such that $B^m$ is spanned. Without loss of generality, we may assume $m = 1$. Suppose that $\{s_0, s_1, \cdots, s_N\}$ is a basis of $H^0(X, B)$. The dimensions of the fibres of the holomorphic map $\Phi_B(x) = [s_0(x) : s_1(x) : \cdots, s_N(x)]$ are not more than $k$. For any sufficiently small open subset $U$ of $X$, let $\theta : B|_U \xrightarrow{\cong} U \times \mathbb{C}$ be a trivialization. We can define a hermitian metric $H$ on $B$ by

$$H(v, v) = \frac{e^{-f(x)}|\theta(v)|^2}{|\theta(s_0(x))|^2 + \cdots + |\theta(s_N(x))|^2}, \quad for \ v \in B|_x.$$

Then $H$ is a smooth hermitian metric. The curvature of $H$ is

$$\begin{aligned} \Theta_H(B) &= -\partial\bar{\partial}\log H \\ &= \partial\bar{\partial}f + \partial\bar{\partial}\log(|\theta(s_0(x))|^2 + \cdots + |\theta(s_N(x))|^2) \\ &= \partial\bar{\partial}f + \Phi_L^* \omega_{FS}, \end{aligned}$$

where $\omega_{FS}$ is the Fubini-Study metric of $\mathbb{P}^N$.

Choose a function $f$ on $X$ as in Lemma 4.2 for the holomorphic map $\Phi_B$ and the Fubini-Study metric $\omega_{FS}$ on $\mathbb{P}^N$. Then $H$ is a smooth hermitian metric such $\Theta_H(B)$ is semipositive and has at least $n - k$ positive eigenvalues. Therefore $B$ is $k$-positive. □

**Corollary 4.4.** *If $E$ is a $k$-ample vector holomorphic bundle of arbitrary rank on a compact complex manifold $X$, then it is $k$-positive.*

**Corollary 4.5.** *Let $E$ be a $k$-ample holomorphic vector bundle on a compact Kähler manifold $X$, then $H^{p,q}(X, E) = 0$ for $p + q \geq n + r + k$.*

*Proof.* By definition $E$ is $k$-ample if and only if the associated tautological line bundle $\xi_E$ is $k$-ample. We know that $\mathbb{P}(E)$ is a Kähler manifold when $X$ is Kähler. Thus Corollary 4.5 follows from Corollary 4.4 together with the Gigante-Girbau vanishing theorem by using Griffiths-Le Potier-Schneider isomorphism in the last part of the proof of Theorem 3.9 (See also Remark 5.9). □

# 5 Vanishing theorems for tensor powers of a Griffiths $k$-positive vector bundles

In this section we will give some vanishing theorems for tensor powers of the Griffiths $k$-positive vector bundles. Similar results for ample vector bundles were obtained by Demailly in [8], nevertheless, using purely algebraic methods. Firstly we will fix some notions for flag manifolds and flag bundles. Our presentation for the geometry of flag bundles is a supplement of those presented by Demailly in [8]. After that we will generalize Theorem 3.8 to the case where the tautological quotient bundles are on arbitrary flag bundles instead of projective bundles. Then we recall Demailly's generalizations of the Griffiths-Le Potier-Schneider isomorphism. Finally we will generalize the vanishing theorems for tensor powers that established by Demailly in [8].

Let $V$ be a complex vector space of dimension $r$. Given any sequence of integers $s = (s_0, s_1, \cdots, s_m)$ such that $0 = s_0 < s_1 < \cdots < s_m = r$, we may consider the manifold $\mathscr{F}_s(V)$ of flags

$$V = V_{s_0} \supset V_{s_1} \supset \cdots \supset V_{s_m} = \{0\}, \quad \text{codim}_\mathbb{C} V_{s_j} = s_j.$$

If $s_j = j$ for all $1 \leq j \leq r$, $\mathscr{F}_s(V)$ is denoted by $\mathscr{F}(V)$, called a complete flag manifold; otherwise called an incomplete flag manifold. The trivial bundle $\mathscr{F}_s(V) \times V$ is denoted by $\widetilde{V}$ and it has flags of tautological subbundles

$$\widetilde{V} = \widetilde{V}_{s_0} \supset \widetilde{V}_{s_1} \supset \cdots \supset \widetilde{V}_{s_m} = \{0\}, \quad \text{rank}_\mathbb{C} \widetilde{V}_{s_j} = r - s_j.$$

Let $Q_{s_j} = \widetilde{V}/\widetilde{V}_{s_j}$ be the tautological quotient vector bundles and denote $P_{s_j} = \widetilde{V}_{s_{j-1}}/\widetilde{V}_{s_j}$ for $1 \leq j \leq m$. Then the cotangent bundle of $\mathscr{F}_s(V)$ is identified by

$$T^*\mathscr{F}_s(V) = \oplus_{j<k}^m \mathrm{Hom}(P_{s_j}, P_{s_k}) \cong \oplus_{j<k}^m P_{s_j}^* \otimes P_{s_k}.$$

So

$$\Omega^p_{\mathscr{F}_s(V)} = \oplus_{p_{1,2}+\cdots+p_{m-1,m}=p} \wedge^{p_{1,2}}(P_{s_1}^* \otimes P_{s_2}) \wedge \cdots \wedge^{p_{m-1,m}}(P_{s_{m-1}}^* \otimes P_{s_m})$$

Using $\Gamma^\rho V^* = \Gamma^{\chi(\rho)} V$ and the Cauchy formula ([27], Sect. 2.1; [28], Sect. 2.2.1) we have

$$\wedge^{p_{j,k}}(P_{s_j}^* \otimes P_{s_k}) = \oplus_{|\rho_{j,k}|=p_{j,k}} \Gamma^{\rho_{j,k}} P_{s_j}^* \otimes \Gamma^{\widetilde{\rho}_{j,k}} P_{s_k} = \oplus_{|\rho_{j,k}|=p_{j,k}} \Gamma^{\chi(\rho_{j,k})} P_{s_j} \otimes \Gamma^{\widetilde{\rho}_{j,k}} P_{s_k}.$$

we then get

$$\Omega^p_{\mathscr{F}_s(V)} = \oplus_{\sum p_{j,k}=p} \otimes_{|\rho_{j,k}|=p_{j,k}}^{1 \leq j \leq m} \Gamma^{\chi(\rho_{1,j})} P_{s_j} \otimes \cdots \otimes \Gamma^{\chi(\rho_{j-1,j})} P_{s_j} \otimes \Gamma^{\widetilde{\rho}_{j,j+1}} P_{s_j} \otimes \cdots \Gamma^{\widetilde{\rho}_{j,m}} P_{s_j}.$$

Since the rank of $P_{s_j}$ and $P_{s_k}$ are $s_j - s_{j-1}$ and $s_k - s_{k-1}$ respectively, we know that the Young diagram of $\rho_{j,k}$ has at most $s_j - s_{j-1}$ rows and $s_k - s_{k-1}$ columns (we refer the reader to ([28], Sect. 2.2.2) for notions used here). In particular we have $\det(P_{s_j} \otimes P_{s_k}^*) = (\det P_{s_j})^{s_k-s_{k-1}} \otimes (\det P_{s_k})^{s_{j-1}-s_j}$. Furthermore $\sum_{j<k} |\rho_{j,k}| \leq \sum_{i<j}(s_j - s_{j-1})(s_k - s_{k-1}) = N_s := \dim \mathscr{F}_s(V)$. Thus for $p = N_s$ we get

$$K_{\mathscr{F}_s(V)} = \otimes_{1 \leq j < k \leq m} \det(P_{s_j} \otimes P_{s_k}^*) = \otimes_{j=1}^m (\det P_{s_j})^{s_{j-1}+s_j-r} := P_s^{c_s},$$

where $c_s = (s_1 - r, \cdots, s_{j-1}+s_j-r, \cdots, s_{m-1}) \in \mathbb{Z}^m$. Consider the flag manifold $\mathscr{F}_s(V)$ as a $Gl(V)$-homogeneous space and write it as $Gl(V)/B_s$, where $B_s$ is a parabolic subgroup of $G$. Then the action of $Gl(V)$ on the tangential bundle $T\mathscr{F}_s(V)$ induces the representation of $Gl(V)$ on $\Omega^p_{\mathscr{F}_s(V)}$ and hence on the canonical line bundle $K_{\mathscr{F}_s(V)}$. Denote $\mathfrak{gl}(V)$ and $\mathfrak{b}_s$ the Lie algebra of $Gl(V)$ and $B_s$ respectively. Then the fibres of $\Omega^p_{\mathscr{F}_s(V)}$ and $K_{\mathscr{F}_s(V)}$ over any point in $\mathscr{F}_s(V)$ could be identified with $\wedge^p(\mathfrak{gl}(V)/\mathfrak{b}_s)^*$ and $\wedge^{N_s}(\mathfrak{gl}(V)/\mathfrak{b}_s)^*$ respectively. In the case of complete flag, $B_s$ is denoted by $B$ and is a Borel subgroup of $Gl(V)$; its Lie algebra is denoted by $\mathfrak{b}$. Now $c_s$ is just the weight associated to the co-adjoint representation on $\wedge^{N_s}(\mathfrak{gl}(V)/\mathfrak{b}_s)^*$, called it canonical weight of $\mathscr{F}_s(V)$. Note that $\det P_{s_j} = \frac{\det Q_{s_j}}{\det Q_{s_{j-1}}}$ we get

$$K_{\mathscr{F}_s(V)} = (\det Q_{s_1})^{s_0-s_2} \otimes \cdots \otimes (\det Q_{s_{m-1}})^{s_{m-2}-s_m} \otimes (\det Q_{s_m})^{s_{m-1}}. \tag{5.1}$$

Since $\det Q_{s_m}$ is a trivial bundle and the Neron-Severi group $\mathrm{NS}(\mathscr{F}_s(V)) \cong \mathbb{Z}^{\oplus m-1}$, we also write the canonical bundle as

$$K_{\mathscr{F}_s(V)} = \mathcal{O}_{\mathscr{F}_s(V)}(s_0 - s_2, \cdots, s_{m-2} - s_m).$$

In particular, in the case of the Grassmannian bundle $\mathrm{Gr}(V,d)$ we have $s = (s_0, s_1, s_2) = (0, d, n)$ and hence $K_{\mathrm{Gr}(V,d)} = (\det Q_{s_1})^{-n}$ by (5.1), usually one denoted it by $\mathcal{O}(-n)$. In particular, for a projective bundle we have $K_{\mathbb{P}(V)} = \mathcal{O}(-n)$, a symbol used popularly in the literatures. This is the main reason why we use several notations different from those used by Demailly in studying the geometry of flag manifolds.

In the case of the complete flag manifold, all $P_j$ are line bundles, we have $T^*\mathscr{F}(V) = \oplus_{j<k}^m (P_j^{-1} \otimes P_k)$ and

$$\Omega^p_{\mathscr{F}(V)} = \wedge^p(\oplus_{j<k}^m P_j^{-1} \otimes P_k) = \oplus_{u \in \mathbb{Z}^r} \nu(u,p) P^u, \tag{5.2}$$

where the weight $u$ and the multiplicity $\nu(u,p)$ are those of $\wedge^p(\mathfrak{gl}(V)/\mathfrak{b})$. In particular,

$$K_{\mathscr{F}(V)} = \otimes_{j=1}^r P_j^{2j-r-1} := P^{c_w}, \tag{5.3}$$

where $c_w = 2(1,2,\cdots,r) - (r+1)(1,1,\cdots,1) = (1-r, 3-r, \cdots, r-1)$.

**Theorem 5.1.** (Bott, [4]) *For any $a = (a_1, a_2, \cdots, a_r) \in \mathbb{Z}^r$ define $\hat{a} = (a - \frac{1}{2}c_w)^{\geq} + \frac{1}{2}c_w$, where $(a - \frac{1}{2}c_w)^{\geq}$ is the sequence obtained by arranging the terms of $(a - \frac{1}{2}c_w)$ in weakly decreasing order. Denote $n(a)$ the number of strict inversions of the order. Then*

$$H^q(\mathscr{F}(V), P^a) = \delta_{q,n(a)} \Gamma^{\hat{a}} V. \tag{5.4}$$

For $a_s = (a_{s_1}, \cdots, a_{s_m}) \in \mathbb{Z}^m$, we define $a \in \mathbb{Z}^r$ such that $a_{s_{j-1}+1} = \cdots = a_{s_j} := a_{s_j}$ for $1 \leq j \leq m \leq r$. For the incomplete flag manifold $\mathscr{F}_s(V)$, there is a natural projection $\eta : \mathscr{F}(V) \to \mathscr{F}_s(V)$ whose fibre is a product of complete flag manifolds $\mathscr{F}(V/V_{s_1}) \times \mathscr{F}(V_{s_1}/V_{s_2}) \times \cdots \mathscr{F}(V_{s_{m-1}}/V_{s_m})$, and $\eta^* P_{s_j} = P_{s_{j-1}+1} \otimes \cdots \otimes P_{s_j}$ and $\eta^* P_s^{a_s} = P^a$.

**Corollary 5.2.** *For the incomplete flag manifold $\mathscr{F}_s(V)$,*

$$H^q(\mathscr{F}_s(V), P_s^{a_s}) = \delta_{q,n(a)} \Gamma^{\hat{a}} V. \tag{5.5}$$

Let $\mathrm{Gr}(V,d)$ denote the Grassmannian of subspaces of $V$ of codimension $d$. In particular $\mathrm{Gr}(V,1) = \mathbb{P}(V)$. If $W$ is a subspace of $V$ of codimension $d$, then the kernel of the map from $\wedge^d V$ to $\wedge^d(V/W)$ is a codimension 1 subspace of $\wedge^d V$. Assigning this codimension 1 subspace of $\wedge^d V$ gives a map $\mathrm{Gr}(V,d) \to \mathbb{P}(\wedge^d E)$, is called the Plücker embedding. An element in $\wedge^d(V/W)$ can be given by a non-zero decomposable $d$-vector

$$\Lambda_{V/W} = v^1 \wedge v^2 \wedge \cdots \wedge v^d,$$

defined up to a constant factor.

The flag manifold $\mathscr{F}_s(V)$ of flags $V = V_{s_0} \supset V_{s_1} \supset \cdots \supset V_{s_m} = \{0\}$, with $\mathrm{codim}_\mathbb{C} V_{s_i} = s_i$, considered as a submanifold of $\mathrm{Gr}(V,s_1) \times \cdots \times \mathrm{Gr}(V,s_m)$, via the Plücker embedding, is an embedded submanifold of the product projective spaces $\mathbb{P}(\wedge^{s_1} V) \times \cdots \times \mathbb{P}(\wedge^{s_m} V)$, we denote this embedding map by $\phi$. Let $\pi_j$ denote the projection from $\mathbb{P}(\wedge^{s_1} V) \times \cdots \times \mathbb{P}(\wedge^{s_m} V)$ to its $j$-th factor $\mathbb{P}(\wedge^{s_j} V)$, and $\xi_{s_j}$ the tautological line bundle of $\mathbb{P}(\wedge^{s_j} V)$. Then we have

$$\phi^* \pi_j^* \xi_{s_j} = \wedge^{s_j}(\widetilde{V}/\widetilde{V_{s_j}}) = \det(Q_{s_j}).$$

Let $E$ be a holomorphic vector bundle of rank $r$ on a compact complex manifold $X$. For every sequence of integers $0 = s_0 < s_1 < \cdots < s_m = r$, we associate to $E$ its flag bundle $\mathscr{F}_s(E) \xrightarrow{\pi} X$ of flags of subbundles $E_{s_0} \supset \cdots \supset E_{s_m}$ with $\mathrm{corank}_{\mathscr{O}_X} E_{s_j} = s_j$. It has tautological quotient vector bundle $Q_{s_j} = (\pi^* E)/(\pi^* E_{s_j}) = \pi^*(E/E_{s_j})$.

**Theorem 5.3.** *If $E$ is Griffiths $k$-positive hermitian holomorphic vector bundle with hermitian metric $h$. Then $\det(Q_{s_1}) \otimes \cdots \otimes \det(Q_{s_{m-1}})$ is $k$-positive.*

*Proof.* We use similar ideas as in Theorem 3.8. Note that $Q_{s_j}^* \cong \pi^*(E_{s_j}^\perp) \subset \pi^* E^*$. There is a naturally associated hermitian metric $H$ on $Q = \det(Q_{s_1}) \otimes \cdots \otimes \det(Q_{s_{m-1}})$. Choose a basis $e^1, \cdots, e^r$ of $E^*$ such that for the flag $(E/E_{s_1})^* \subset \cdots \subset (E/E_{s_m})^*$, we have $(E/E_{s_j})^* = \mathrm{Span}_{\mathbb{C}}\{e^1, \cdots, e^{s_j}\}$. We use $\Lambda_{E_{s_j}/E_{s_{j+1}}} = e^{s_j+1} \wedge \cdots \wedge e^{s_{j+1}}$ to represent the Plücker coordinate of $(E_{s_j}/E_{s_{j+1}})^*$. For $(z; \Lambda) = (z; [\Lambda_{E_{s_1}/E_{s_2}}, \cdots, \Lambda_{E_{s_{m-1}}/E_{s_m}}]) \in Q_{s_1}^* \otimes \cdots \otimes Q_{s_{m-1}}^*$, the hermitian metric $H$ is defined by

$$H(\Lambda, \Lambda) = \sum_{j=1}^{m-1} (\Lambda_{E_{s_j}/E_{s_{j+1}}}, \Lambda_{E_{s_j}/E_{s_{j+1}}})_z = \sum_{j=1}^{m-1} |\Lambda_{E_{s_j}/E_{s_{j+1}}}|_z^2,$$

where

$$(\Lambda_{E_{s_j}/E_{s_{j+1}}}, \Lambda_{E_{s_j}/E_{s_{j+1}}})_z = \det\left((h(e^k, e^l))_{s_j < k, l \leq s_{j+1}}\right).$$

Choose local normal coordinate at $z$ such that $h(0) = \delta_{jk}$ and $dh(0) = 0$. Then

$$\partial\bar{\partial}(\Lambda_{E_{s_j}/E_{s_{j+1}}}, \Lambda_{E_{s_j}/E_{s_{j+1}}})_z$$
$$= (\mathrm{Tr}_{E_{s_j}/E_{s_{j+1}}}(\Theta_h(E^*)))(\Lambda_{E_{s_j}/E_{s_{j+1}}}, \Lambda_{E_{s_j}/E_{s_{j+1}}})_z + (\partial\Lambda_{E_{s_j}/E_{s_{j+1}}}, \partial\Lambda_{E_{s_j}/E_{s_{j+1}}})_z,$$

where

$$\mathrm{Tr}_{E_{s_j}/E_{s_{j+1}}}(\Theta_h(E^*)) = \sum_{k=s_j+1}^{s_{j+1}} \frac{\Theta_h(E^*)(e^k, e^k)}{|e^k|^2}.$$

Thus

$$i\Theta_H(Q) = \frac{\sum_{j=1}^{m-1}(\mathrm{Tr}_{E_{s_j}/E_{s_{j+1}}}(\Theta_h(E)))|\Lambda_{E_{s_j}/E_{s_{j+1}}}|^2}{\sum_{j=1}^{m-1}|\Lambda_{E_{s_j}/E_{s_{j+1}}}|^2}$$
$$+ \frac{\sum_{j=1}^{m-1}(|\Lambda_{E_{s_j}/E_{s_{j+1}}}|^2|\partial\Lambda_{E_{s_j}/E_{s_{j+1}}}|^2 - (\partial\Lambda_{E_{s_j}/E_{s_{j+1}}}, \Lambda_{E_{s_j}/E_{s_{j+1}}})(\Lambda_{E_{s_j}/E_{s_{j+1}}}, \partial\Lambda_{E_{s_j}/E_{s_{j+1}}}))}{(\sum_{j=1}^{m-1}|\Lambda_{E_{s_j}/E_{s_{j+1}}}|^2)^2}.$$

The second term of right hand side is equal to $\sum_{j=1}^{m-1} \kappa_j \theta_j$ with $0 < \kappa_j < 1$, and

$$\theta_j = \frac{|\Lambda_{E_{s_j}/E_{s_{j+1}}}|^2 |\partial\Lambda_{E_{s_j}/E_{s_{j+1}}}|^2 - (\partial\Lambda_{E_{s_j}/E_{s_{j+1}}}, \Lambda_{E_{s_j}/E_{s_{j+1}}})(\Lambda_{E_{s_j}/E_{s_{j+1}}}, \partial\Lambda_{E_{s_j}/E_{s_{j+1}}})}{|\Lambda_{E_{s_j}/E_{s_{j+1}}}|^4}$$

is the curvature of tautological quotient line bundle $\det(E_{s_j}/E_{s_{j+1}})$ of the Grassmannian manifold $\mathrm{Gr}(E_{s_j}, s_{j+1} - s_j)$ (cf. [6]), therefore $\theta_j$ is positive of rank $s_j(s_{j+1} - s_j)$. Since the dimension of the fibre of the flag bundle $\mathscr{F}_s(E)$ is $N_s = s_1(s_2 - s_1) + s_2(s_3 - s_2) + \cdots + s_{m-1}(r - s_{m-1})$, it follows that $\sum_{j=1}^{m-1} \kappa_j \theta_j$ is positive definite along the fibres. Hence $\det(Q_{s_1}) \otimes \cdots \otimes \det(Q_{s_{m-1}})$ is $k$-positive if $E$ is Griffiths $k$-positive. $\square$

**Corollary 5.4.** *Let $E$ be a Griffiths k-positive hermitian holomorphic vector bundle. Then*

*(i). The line bundle $Q_s^{a_s} = \det(Q_{s_1})^{a_1} \otimes \cdots \otimes \det(Q_{s_{m-1}})^{a_{m-1}}$, also denoted by $\mathcal{O}_{\mathscr{F}_s(E)}(a_1, \cdots, a_{m-1})$, is semi-positive if $a_j \geq 0$ and k-positive if $a_j > 0$ for any $1 \leq j \leq m-1$;*

*(ii). The line bundle $P_s^{a_s} = \det(P_{s_1})^{a_1} \otimes \cdots \otimes \det(P_{s_{m-1}})^{a_{m-1}}$ is semi-positive if $a_j \geq a_{j+1}$ and $a_{m-1} \geq 0$; k-positive if $a_j > a_{j+1}$ and $a_{m-1} > 0$, for any $1 \leq j \leq m-2$.*

*Proof.* (i) is proved similarly as Theorem 5.3. For (ii), note that

$$\begin{aligned} P_s^{a_s} &= \det(Q_{s_1})^{a_1-a_2} \otimes \cdots \otimes \det(Q_{s_{m-2}})^{a_{m-2}-a_{m-1}} \otimes \det(Q_{s_{m-1}})^{a_{m-1}} \\ &= \mathcal{O}_{\mathscr{F}_s(E)}(a_1-a_2, \cdots, a_{m-2}-a_{m-1}, a_{m-1}). \end{aligned}$$

□

In particular, if the base manifold $X$ is shrunk to a point we get

**Corollary 5.5.** *On the flag manifold $\mathscr{F}_s(V)$, the line bundle $Q_s^{a_s} = \det(Q_{s_1})^{a_1} \otimes \cdots \otimes \det(Q_{s_m})^{a_m}$ is ample if $a_j > 0$ for any $1 \leq j \leq m-1$; and the line bundle $P_s^{a_s} = \det(P_{s_1})^{a_1} \otimes \cdots \otimes \det(P_{s_m})^{a_m}$ is ample if $a_j > a_{j+1}$ for any $1 \leq j \leq m-1$.*

Denote
$$\begin{aligned} \eta^* \Omega^p_{\mathscr{F}_s(V)} &= \oplus_u \nu_s(u,p) P^u, \\ \eta^* (\Omega^p_{\mathscr{F}_s(V)})^* &= \oplus_u \nu_s(u',p) P^{u'}, \end{aligned} \tag{5.6}$$

where the weight $u$ (resp. $u'$) and multiplicity $\nu(u,p)$ (resp. $\nu(u',p)$) are those of $\wedge^p(\mathfrak{gl}(V)/\mathfrak{b}_s)^*$ (resp. $\wedge^p(\mathfrak{gl}(V)/\mathfrak{b}_s)$).

**Lemma 5.6.** *(cf. [8], Lemma 2.21). The weight $u$ of $\wedge^p(\mathfrak{gl}(V)/\mathfrak{b}_s)^*$ verifies*

$$u_\mu - u_\lambda \leq \min\{p+1, r+1-(\mu-\lambda), N_s - p + (s_{j+1}-s_{j-1})\}$$

*for $s_{j-1} < \lambda \leq s_j < \mu \leq s_{j+1}$, for $1 \leq j \leq m-1$.*

Using Corollary 5.5 and Lemma 5.6, it is easy to prove

**Lemma 5.7.** *(cf. [8], Lemma 3.4). Assume the weight $a$ satisfying*
$$\begin{aligned} &a_{s_j} - a_{s_{j+1}} \geq 1 \quad \text{if} \quad p = N_s \quad \text{and otherwise} \\ &a_{s_j} - a_{s_{j+1}} \geq \min\{p, N_s - p + (s_{j+1}-s_j) - 1, r+1-(s_{j+1}-s_{j-1})\}, \end{aligned} \tag{5.7}$$
*then the direct image sheaf*

$$R^q \pi_* (\Omega^p_{\mathscr{F}_s(E)/X} \otimes P_s^{a_s}) = 0, \quad for \quad q \geq 1 \tag{5.8}$$

$$\pi_*(\Omega^p_{\mathscr{F}_s(E)/X} \otimes P_s^{a_s}) = \oplus_u \nu_s(u,p) \Gamma^{a+u} E. \tag{5.9}$$

*In particular we have*

$$\pi_*(P_s^{a_s}) = \Gamma^a E, \qquad \pi_*(\Omega^{N_s}_{\mathscr{F}_s(E)/X} \otimes P_s^{a_s}) = \Gamma^{a+c_w} E.$$

The Leray spectral sequence implies therefore the following Griffiths-Le Potier-Schneider-Demailly isomorphism:

**Theorem 5.8.** (cf. [8], Theorem 3.8). *Under assumption (5.7), for any holomorphic line bundle $B$ on $X$, we have for all $q \geq 0$:*

$$H^q(\mathscr{F}_s(E), G^{p,p+t} \otimes P_s^{a_s} \otimes \pi^* B) \cong \oplus_u \nu_s(u,t) H^q(X, \Omega_X^p(\Gamma^{a+u} E \otimes B)), \tag{5.10}$$

*where $G^{p,p+t} = \pi^*(\Omega_X^p) \otimes \Omega_{\mathscr{F}_s(E)/X}^t$.*

In particular, if $t = N_s$, then

$$H^q(\mathscr{F}_s(E), G^{p,p+N_s} \otimes P_s^{a_s} \otimes \pi^* B) \cong H^q(X, \Omega_X^p(\Gamma^{a+c_w} E \otimes B)) \tag{5.11}$$

when $a_{s_j} - a_{s_{j+1}} \geq 1$ and $q \geq 0$.

It is not difficult to prove that the canonical bundle of the flag bundle $\mathscr{F}_s(E)$ and that of $X$ are related by

$$K_{\mathscr{F}_s(E)} = \pi^*(K_X \otimes (\det E)^{s_{m-1}}) \otimes \mathcal{O}_{\mathscr{F}_s(E)}(s_0 - s_2, s_1 - s_3, \cdots, s_{m-2} - s_m).$$

Therefore we have the following generalization of the Griffiths isomorphism

$$H^q(\mathscr{F}_s(E), K_{\mathscr{F}_s(E)} \otimes P_s^{a_s} \otimes \pi^* B) \cong H^q(X, K_X \otimes \Gamma^{a+c_w} E \otimes B). \tag{5.12}$$

**Remark 5.9.** In (5.10) if we take $t = p$ then we get the isomorphism for Dolbeault cohomology groups

$$H^q(Y, \Omega_Y^p(P_s^{a_s} \otimes \pi^* B)) \cong \oplus_u \nu_s(u,p) H^q(X, \Omega_X^p(\Gamma^{a+u} E \otimes B)). \tag{5.13}$$

For projective bundles, the following two special cases of the Griffiths-Le Potier-Schneider-Demailly isomorphism are frequently used.

The first case is $p = 0$ and $a_s = (l, 0)$, the later corresponds to $a = (l, 0, \cdots, 0)$. We have $P_s^{a_s} = \xi_E^l$. Therefore

$$H^q(\mathbb{P}(E), \xi_E^l \otimes \pi^* B)) \cong H^q(X, S^l E \otimes B)),$$

which was proved by Griffiths in [17]. This together with the formula $K_{\mathbb{P}(E)} = \xi_E^{-r} \otimes \pi^*(\det E \otimes K_X)$ give

$$H^q(\mathbb{P}(E), K_{\mathbb{P}(E)} \otimes \xi_E^l \otimes \pi^* B)) \cong H^q(X, K_X \otimes S^{l-r} E \otimes \det E \otimes B))$$

if $l \geq r$. It is a special case of (5.12) for $a_s = (l, 0)$. Since $c_s = (1 - r, 1)$, we have $a_s + c_s = (l - r + 1, 1)$, and accordingly $a + c = (l - r + 1, 1, \cdots, 1)$ and $\Gamma^{a+c} E = S^{l-r} E \otimes \det E$. Taking $l = r + 1$ in the isomorphism above and using the Gigante-Girbau vanishing theorem together with Proposition 3.7, we get another proof of the case a) of Corollary 3.13.

The second special case is $a_s = (1, 0)$. Then there is no weight $u$ for an irreducible representation such that $a + u$ is non-increasing. Therefore for any $1 \leq p, q \leq n$, we have

$$H^q(\mathbb{P}(E), \Omega^p_{\mathbb{P}(E)}(\xi_E \otimes \pi^* B)) \cong H^q(X, \Omega^p_X(E \otimes B)),$$

which was firstly proved by Le Potier [24, 19] and a simplified proof was given by Schneider [37, 38].

**Theorem 5.10.** *Let $E$ be a Griffiths $k$-positive holomorphic vector bundle of rank $r$ and $B$ a semipositive holomorphic line bundle on a compact Kähler manifold $X$. Let $\Gamma^a E$ be an irreducible tensor power representation of $Gl(E)$ of highest weight $a \in \mathbb{Z}^r$ with $h \in \{1, \cdots, r-1\}$ such that $a_1 \geq a_2 \geq \cdots \geq a_h > a_{h+1} = \cdots = a_r = 0$. Then*

$$H^{n,q}(X, \Gamma^a E \otimes (\det E)^h \otimes B) = 0 \quad for \ q > k. \tag{5.14}$$

*Proof.* Define $s_1 < s_2 < \cdots < s_{m-1}$ as the sequence of indices $\lambda \in \{1, \cdots, r-1\}$ such that $a_\lambda > a_{\lambda+1}$. Then $s_{m-1} = h$ and $s_m = r$. Let $a_s = (a_{s_1}, \cdots, a_{s_m})$ and $a'_s = a_s + (h, \cdots, h) - c_s$. Recall that $c_s = (s_1 - r, \cdots, s_{j-1} + s_j - r, \cdots, s_{m-1})$, therefore the canonical weight $c_s$ is non-decreasing and $a'_{s_1} > a'_{s_2} > \cdots > a'_{s_m} = 0$. Thus $P_s^{a'_s}$ is $k$-positive by Corollary 5.4. Note $\Gamma^{a'_s + c_s} E = \Gamma^a E \otimes (\det E)^h$. By (5.12) we have

$$H^{n,q}(X, \Gamma^a E \otimes (\det E)^h \otimes B) \cong H^{n+N_s,q}(\mathscr{F}_s(E), P_s^{a'_s} \otimes \pi^* B).$$

Note that $\dim \mathscr{F}_s(E) = n + N_s$ and since $\pi^* B$ is semipositive and $P_s^{a'_s}$ is $k$-positive we know $P_s^{a'_s} \otimes \pi^* B$ is $k$-positive by Proposition 3.6 b), thus the cohomology group of the right hand side is zero by the Gigante-Girbau vanishing theorem. $\square$

**Theorem 5.11.** *Let $E$ and $B$ be the vector bundles as in Theorem 5.8. Then for all integers $p + q > n + k$, $l \geq 1$, $m \geq n - p + r - 1$, we have*

$$H^{p,q}(X, E^{\otimes l} \otimes (\det E)^m \otimes B) = 0. \tag{5.15}$$

*Proof.* Since $E^{\otimes l}$ is decomposed into a direct sum of irreducible representations $\Gamma^a E$ of $Gl(E)$, it suffices to prove that

$$H^{p,q}(X, \Gamma^a E \otimes (\det E)^m \otimes B) = 0 \tag{5.16}$$

for $p + q > n + k, m \geq n - p + r - 1$ and $a_1 \geq \cdots \geq a_r = 0$ not all zero. We prove it by backward induction on $p$. The case $p = n$ is already established by Theorem 5.10. Define $s$ as in the proof of Theorem 5.10, and set $a'_s = a_s + (m, \cdots, m) - c_s$. Then $P_s^{a'_s}$ is $k$-positive. By (5.11) we get

$$H^{p,q}(X, \Gamma^a E \otimes (\det E)^m \otimes B) \cong H^q(Y, G^{p,p+N_s} \otimes P_s^{a'_s} \otimes \pi^* B). \tag{5.17}$$

Put $F^{p,t} = \pi^*(\Omega_X^p) \wedge \Omega_{\mathscr{F}_s(E)}^{t-p}$. Then we have the following exact sequence

$$0 \to F^{p+1,t} \to F^{p,t} \to G^{p,t} \to 0. \tag{5.18}$$

Note that $F^{p,p+N_s} = \Omega_{\mathscr{F}_s(E)}^{p+N_s}$, so we have

$$0 \to F^{p+1,p+N_s} \to \Omega_{\mathscr{F}_s(E)}^{p+N_s} \to G^{p,p+N_s} \to 0. \tag{5.19}$$

The Gigante-Girbau vanishing theorem applied to $P_s^{a'_s} \otimes \pi^* B$ yields

$$H^q(Y, \Omega_{\mathscr{F}_s(E)}^{p+N_s} \otimes P_s^{a'_s} \otimes \pi^* B) = 0, \quad for\ p+q > n+k.$$

The cohomology group in (5.17) will therefore vanishes if and only if

$$H^{q+1}(\mathscr{F}_s(E), F^{p+1,p+N_s} \otimes P_s^{a'_s} \otimes \pi^* B) = 0. \tag{5.20}$$

By (5.18) we have $F^{p,t}/F^{p+1,t} = G^{p,t}$ and hence $F^{p+1,p+N_s}$ has a filtration with associated graded bundles $\oplus_{t \geq 1} G^{p+t,p+N_s}$. To prove (5.20), it is thus enough to verify

$$H^{q+1}(\mathscr{F}_s(E), G^{p+t,p+N_s} \otimes P_s^{a'_s} \otimes \pi^* B) = 0. \ t \geq 1. \tag{5.21}$$

Note that the direct image sheaf

$$R^l \pi_*(G^{p+t,p+N_s} \otimes P_s^{a'_s} \otimes \pi^* B) = \Omega_X^{p+t} \otimes B \otimes R^l \pi_*(\Omega_{\mathscr{F}_s(E)/X}^{N_s-t} \otimes P^{a_s+(m,\cdots,m)-c_s}).$$

Corollary 5.5 and Bott's theorem yield

$$\begin{aligned} & R^l \pi_*(\Omega_{\mathscr{F}_s(E)/X}^{N_s-t} \otimes P^{a_s+(m,\cdots,m)-c_s})|_{p \in X} \\ = & H^l(\mathscr{F}_s(V), \Omega_{\mathscr{F}_s(V)}^{N_s-t}(P^{a_s+(m,\cdots,m)-c_s})) = \begin{cases} 0 & \text{for } l > t, \text{ otherwise} \\ \oplus_b \Gamma^b E, & b_r \geq m-t, \end{cases} \end{aligned}$$

$$R^l \pi_*(G^{p+t,p+N_s} \otimes P_s^{a'_s} \otimes \pi^* B) = \begin{cases} 0 & \text{for } l > t, \text{ otherwise} \\ \oplus_{b,j} \Omega_X^{p+t} \otimes \Gamma^b E \otimes (\det E)^j \otimes B, \end{cases}$$

where the last sum runs over weights $b$ such that $b_r = 0$ and integer $j$ such that $j \geq m - t \geq n - (p+t) + r - 1$. Using the Leray spectral sequence, we have

$$H^{q+1}(\mathscr{F}_s(E), G^{p+t,p+N_s} \otimes P_s^{a'_s} \otimes \pi^* B) = \oplus_{l,b,j} H^{p+t,q+1-l}(X, \Gamma^b E \otimes (\det E)^j \otimes B).$$

Note in the second sum, $l \leq t$ and $t \geq 1$, hence $p+t > p$ and $(p+t)+(q+1-l) \geq p+q+1 > n+k$. From the induction hypothesis we have $H^{p+t,q+1-l}(X, \Gamma^b E \otimes (\det E)^j \otimes B) = 0$ and (5.21) is proved. □

**Remark 5.12.** Let $B$ be a semipositive line bundle and $E$ a Griffiths $k$-positive holomorphic vector bundle of rank $r$ on a compact Kähler manifold $X$, and $s = (s_0, s_1, \cdots, s_m)$ an increasing sequence of integers with $s_m = r$. Let $a_s = (a_1, \cdots, a_m)$ be a decreasing sequence of integers satisfying (5.7). Then by (5.13) we have

$$H^q(\mathscr{F}_s(E), \Omega^p_{\mathscr{F}_s(E)}(P_s^{a_s} \otimes \pi^* B)) \cong \oplus_u \nu_s(u, p) H^q(X, \Omega^p_X(\Gamma^{a+u} E \otimes B)).$$

By Corollary 5.4 and the Gigante-Girbau vanishing theorem the cohomology group of left hand side vanishes when $p + q > n + k + N_s$. Thus

$$H^q(X, \Omega^p_X(\Gamma^{a+u} E \otimes B)) = 0, \qquad for \ \ p + q > n + k + N_s, \tag{5.22}$$

where $u$ is any weight of the adjoint representation of $Gl(V)$ on $\wedge^p(\mathfrak{gl}(V)/\mathfrak{b}_s)^*$. This generalizes the last statement of Theorem 3.8. However, its relations to vanishing phenomena in Theorems 5.11 is delicate, we think that Theorem 5.11 includes it only in some particular cases.

# 6 Closing Remarks

As suggested in Remark 2.5, we have the following first question
**Question 6.1.** *Give examples of pseudo-effective line bundles on compact Kähler manifolds, such that the Nadel type vanishing theorem for Dolbeault cohomology are not true.*

Our second question is concerned with positivity of tensor powers of vector bundles. As proved by Proposition 3.6 c), if $E$ is $(k, s)$-positive, then so is $\Gamma^a E \otimes (\det E)^m$. In particular, $S^p E$ and $\wedge^q E$ are $(k, s)$-positive. However, as indicated by Corollary 3.12, Theorem 5.10 and Theorem 5.11, $\Gamma^a E \otimes (\det E)^m$ holds stronger positivity than $(k, s)$-positivity. We note that Corollary 3.13 is based on the critical Proposition 3.10, which is directly established by curvature calculations. Nevertheless, Theorem 5.10 and Theorem 5.11 mainly use the generalized Griffiths-Le Potier-Schneider-Demailly isomorphism. The later method is mysterious and we think that it only recovers partial positivity. So we present the question
**Question 6.2.** *Let $E$ be a $(k, s)$-positive hermitian holomorphic vector bundle on a compact complex manifold. Find out the minimum $p$ the maximum $q$ directly via curvature calculations such that $\Gamma^a E \otimes (\det E)^m$ is $(p, q)$-positive.*

Finally we have the following conjecture as suggested by a combination of Theorem 2.4 and Theorem 3.9.
**Conjecture 6.3.** *Let $X$ be a compact Kähler manifold of dimension $n$ and $E$ a hermitian holomorphic vector bundle of rank $r$ on $X$ such that $E >_{(k,s)} 0$. Then*

$$H^q(X, \Omega^p_X \otimes E) = 0, \quad for \ \ p + q > n + k \ \ and \ \ s \geq \min\{n - q + 1, r\}.$$

Note that Theorem 2.4 and Theorem 3.9 are the special cases of $r = 1$ and $p = n$ respectively. In particular, if Conjecture 6.3 is true, it includes the Kodaira vanishing theorem, the Nakano vanishing theorem, the Gigante-Girbau vanishing theorem and the Demailly vanishing theorem (Theorem 3.9 when $k = 0$) as special cases. However the method of the proof of Theorem 2.4 is incapable of proving Conjecture 6.3 if the rank $r > 1$, the main difficult in using such a method is that we can diagonalize simultaneously the Kähler metric and the curvature $i\Theta_h(E)$ only along one direction in $E$ and it is not enough if $E$ is not a line bundle.

# Acknowledgments

I would like to thank T. Ohsawa for helpful discussions and generous supports, and J.P. Demailly, A. Fujiki, J. Le Potier, A. J. Sommese, S. Takayama and H. Tsuji for answering my questions. This work was done in 2007, during my stay at Nagoya University, I would like to thank the members of Department of Mathematics for their hospitality. Last but not least, I'd like to thank the anonymous referee, for his (her) patience in reading this article, and made very good suggestions.

This work was partly supported by the Fundamental Research Funds for the Central Universities (No. 09lgpy49) and NSFC (No. 10801082).